\documentclass[12pt]{article}%
\usepackage{amsmath,amssymb,fullpage}
\usepackage{amsfonts}
\usepackage{mathrsfs}
\usepackage{amsthm}
\usepackage[ bookmarks=true,         bookmarksnumbered=true, colorlinks=true,
pdfstartview=FitV, linkcolor=blue, citecolor=blue, urlcolor=blue]{hyperref}
\usepackage{amsmath}
\usepackage{amssymb}
\usepackage{graphicx}%
\providecommand{\U}[1]{\protect\rule{.1in}{.1in}}
\newtheorem{example}{Example}[section]
\newtheorem{theorem}[example]{Theorem}
\newtheorem{definition}[example]{Definition}

\newtheorem{proposition}[example]{Proposition}
\newtheorem{problem}[example]{Problem}
\newtheorem{remark}[example]{Remark}

\numberwithin{equation}{section}

\def\1B{\text{1\!\!I}}

\let\Section=\section
\def\section{\setcounter{equation}{0}\Section}
\begin{document}

\title{Optimal control of forward-backward stochastic Volterra equations}
\author{Nacira Agram$^{1,2,3}$, Bernt {\O }ksendal$^{1,2}$ and Samia Yakhlef$^{3}$ }
\date{15 September 2017}
\maketitle

\footnotetext[1]{Department of Mathematics, University of Oslo, P.O. Box 1053
Blindern, N--0316 Oslo, Norway. Email: \texttt{naciraa@math.uio.no},
oksendal@math.uio.no.}

\footnotetext[2]{This research was carried out with support of the Norwegian
Research Council, within the research project Challenges in Stochastic
Control, Information and Applications (STOCONINF), project number 250768/F20.}
\footnotetext[3]{Department of Mathematics, University of Biskra, Algeria.
Email: \texttt{samiayakhelef@yahoo.fr.}}

\noindent\textbf{Abstract:} We study the problem of optimal control of a
coupled system of forward-backward stochastic Volterra equations. We use
Hida-Malliavin calculus to prove a sufficient and a necessary maximum
principle for the optimal control of such systems. Existence and uniqueness of
backward stochastic Volterra integral equations are proved. As an application
of our methods, we solve a recursive utility optimisation problem in a
financial model with memory.\newline\newline



\section{Introduction}

The purpose of this paper is to establish solution techniques for optimal
control of coupled systems of stochastic Volterra equations. Stochastic
Volterra equations appear in models for dynamic systems with noise and memory.
As a motivating example, consider the following Volterra equation, modelling a
stochastic cash flow
%In addition to assumptions from the previous sections, we want to consider an
%optimal recursive utility problem in the similar manner with \cite{AO2} and
%\cite{OS}. Let us consume then at the rate $c(t)$, the cash flow
$X(t)=X^{c}(t)$ subject to a consumption rate $c(t)$ at time $t$:%

\begin{equation}%
\begin{array}
[c]{c}%
X(t)=\xi(t)+\int_{0}^{t}\left(  \alpha(t,s)-c(s)\right)  X(s)ds+\int_{0}%
^{t}\beta(t,s)X(s)dB(s)\\
\text{ \ \ \ \ \ \ \ }+\int_{0}^{t}\int_{%
%TCIMACRO{\U{211d} }%
%BeginExpansion
\mathbb{R}
%EndExpansion
_{0}}\pi(t,s,e)X(s)\tilde{N}(ds,de),t\in\left[  0,T\right]  ,
\end{array}
\label{eq1.1}%
\end{equation}
where $\xi:[0,T]\rightarrow\mathbb{R}$ and $\alpha$, $\beta:\left[
0,T\right]  ^{2}\rightarrow%
%TCIMACRO{\U{211d} }%
%BeginExpansion
\mathbb{R}
%EndExpansion
$ and $\pi:\left[  0,T\right]  ^{2}\times%
%TCIMACRO{\U{211d} }%
%BeginExpansion
\mathbb{R}
%EndExpansion
_{0}\rightarrow%
%TCIMACRO{\U{211d} }%
%BeginExpansion
\mathbb{R}
%EndExpansion
$ are deterministic functions with $\alpha$, $\beta$ and $\pi$
bounded.\newline Here $B(t)=B(t,\omega)$ and $N(dt,de)=N(dt,de,\omega)$ are a
Brownian motion and an independent Poisson random measure, respectively, on a
complete probability space $\left(  \Omega,\mathcal{F},P\right)  $. The
compensated Poisson random measure $\tilde{N}$ is defined by $\tilde
{N}(dt,de)=N(dt,de)-\nu(de)dt$, where $\nu$ is the L\'{e}vy measure of $N$. We
denote by $\mathbb{F}=\{\mathcal{F}_{t}\}_{t\geq0}$ the right-continuous
complete filtration generated by $B$ and $N$ and we let
\[
\mathbb{G}:=\{\mathcal{G}_{t}\}_{t\geq0}%
\]
be a given right-continuous complete subflitration of $\mathbb{F}$, in the
sense that
\[
\mathcal{G}_{t}\subseteq\mathcal{F}_{t},\text{ for all }t\in\lbrack0,T].
\]
The $\sigma$-algebra $\mathcal{G}_{t}$ represents the information available to
the consumer at time $t$. Let $\mathcal{P}(\mathbb{F})$ be the $\sigma
$-algebra of $\mathbb{F}$-predictable subsets of $\Omega\times%
%TCIMACRO{\U{211d} }%
%BeginExpansion
\mathbb{R}
%EndExpansion
_{+},$ i.e., the $\sigma$-algebra generated by the left continuous
$\mathbb{F}$-adapted processes.\newline
%Moreover, for $k\geq1$, we denote by
%$\mathcal{B}(\mathbb{R}^{k})$ the Borel-$\sigma$-field over $\mathbb{R}^{k}$.
%Throughout this work, we introduce some basic notations and spaces for any
%$\sigma$-algebra ${\mathcal{F}}$ and any filtration $\mathbb{F=}\left(
%\mathcal{F}_{t}\right)  _{t\geq0}.$
The \emph{forward stochastic Volterra integral equation} (FSVIE) \eqref{eq1.1}
can be written in differential form as%
\begin{equation}%
\begin{array}
[c]{l}%
dX(t)=\xi^{\prime}(t)dt+\left(  \alpha(t,t)-c(t)\right)  X(t)dt+\left(
\int_{0}^{t}\frac{\partial\alpha}{\partial t}(t,s)X(s)ds\right)  dt\\
+\beta(t,t)X(t)dB(t)+\left(  \int_{0}^{t}\frac{\partial\beta}{\partial
t}(t,s)X(s)dB(s)\right)  dt\\
+\int_{\mathbb{R}_{0}}\pi(t,t,e)X(t)\tilde{N}(dt,de)\\
+\left(  \int_{0}^{t}\int_{\mathbb{R}_{0}}\frac{\partial\pi}{\partial
t}(t,s,e)X(s)\tilde{N}(ds,de)\right)  dt,t\in\left[  0,T\right]  .
\end{array}
\label{eq1.2}%
\end{equation}
From \eqref{eq1.2} we see that the dynamics of $X(t)$ contains history or
memory terms represented by the $ds$-integrals.\newline Following a suggestion
of Duffie and Epstein \cite{DE} we now model the total utility of the
consumption rate $c(t)$ by a \emph{recursive utility} process $Y(t)=Y^{c}(t)$
defined by the equation%

\begin{equation}
Y(t)=\mathbb{E}\left[  -\int_{t}^{T}\left.  \{\gamma(s)Y(s)+\ln
(c(s)X(s))\}ds\right\vert \mathcal{F}_{t}\right]  ,\quad t\in\left[
0,T\right]  .\label{eq1.3}%
\end{equation}
By the martingale representation theorem we see that there exist processes
$Z(t),K(t,e)$ such that the triple $(Y,Z,K)$ solves the backward stochastic
differential equation (BSDE)
\begin{equation}%
\begin{cases}
dY(t)=-\left[  \gamma(t)Y(t)+\ln(c(t)X(t))\right]  dt+Z(t)dB(t)\\
\text{ \ \ \ \ \ \ \ \ \ \ \ \ \ \ \ \ \ }+\int_{\mathbb{R}_{0}}%
K(t,e)\tilde{N}(dt,de),\quad t\in\left[  0,T\right]  ,\\
Y(T)=0.
\end{cases}
\label{eq1.4}%
\end{equation}
We now consider the \emph{optimal recursive utility problem} to maximise the
total recursive utility of the consumption. In other words, we want to find an
optimal consumption rate $c^{\ast}\in\mathcal{U}_{\mathbb{G}}$ such that
\begin{equation}
\sup_{c\in\mathcal{U}_{\mathbb{G}}}Y^{c}(0)=Y^{c^{\ast}}(0),\label{eq1.5}%
\end{equation}
where $\mathcal{U}_{\mathbb{G}}$ is a given set of admissible $\mathbb{G}%
$-adapted consumption processes.\newline This is a problem of optimal control
of a coupled system consisting of the forward stochastic Volterra equation
\eqref{eq1.1} and the BSDE \eqref{eq1.4}. In the following sections we will
present solution methods for general optimal control for systems of
forward-backward stochastic Volterra equations. Then in the last section we
will apply the methods to solve the optimal recursive utility consumption
problem above.\newline There has been a lot of research activity recently
within stochastic Volterra integral equations (SVIE), both of forward and
backward type. See e.g. \cite{AO}, \cite{L}, \cite{SWY}, \cite{SW},
\cite{WZ1}, \cite{WS}, \cite{WZ}, \cite{WX}, \cite{Y} and \cite{Yo}. Perhaps
the paper closest to our paper is \cite{WS}. However, that paper has a
different approach than our paper, does not have a sufficient maximum
principle and does not deal with jumps and partial information, as we do.

\section{Stochastic maximum principle for FBSVE}

This section is an extension to forward-backward systems of the results
obtained in \cite{AO}. We consider a system governed by a coupled system of
controlled forward-backward stochastic Volterra equations (FBSVE) of the form:%
\begin{align}
&
\begin{array}
[c]{c}%
X(t)=\xi(t)+\int_{0}^{t}b(t,s,X(s),u(s))ds+\int_{0}^{t}\sigma
(t,s,X(s),u(s))dB(s)\\
\text{ \ \ \ \ \ \ \ \ \ \ \ \ \ \ \ \ \ \ }+\int_{0}^{t}\int_{%
%TCIMACRO{\U{211d} }%
%BeginExpansion
\mathbb{R}
%EndExpansion
_{0}}\theta(t,s,X(s),u(s),e)\tilde{N}(ds,de),t\in\left[  0,T\right]  ,
\end{array}
\label{a1}\\
&
\begin{array}
[c]{l}%
Y(t)=\eta(X(T))+\int_{t}^{T}g(t,s,X(s),Y(s),Z(t,s),K(t,s,\cdot),u(s))ds\\
\text{ \ \ \ \ \ \ \ \ \ \ \ \ \ \ \ \ \ \ \ }-\int_{t}^{T}Z(t,s)dB(s)-\int
_{t}^{T}\int_{%
%TCIMACRO{\U{211d} }%
%BeginExpansion
\mathbb{R}
%EndExpansion
_{0}}K(t,s,e)\tilde{N}(ds,de),t\in\left[  0,T\right]  .
\end{array}
\label{a2}%
\end{align}
The quadruple $(X,Y,Z,K)$ is said to be a solution of \eqref{a1}-\eqref{a2} if
it satisfies both equations. To the best of our knowledge, results about
existence and uniqueness of solutions for such general systems are not known.
Conditions under which there exists a unique solution $\left(  Y,Z,K\right)  $
of $\left(  \ref{a2}\right)  $ are studied in section $3$.\newline In the
above, the functions $\xi,\eta$ are assumed to be deterministic and $C^{1}$,
while the functions%
\[%
\begin{array}
[c]{ll}%
b(t,s,x,u) & :\left[  0,T\right]  ^{2}\times%
%TCIMACRO{\U{211d} }%
%BeginExpansion
\mathbb{R}
%EndExpansion
\times\mathbb{U}\times\Omega\rightarrow%
%TCIMACRO{\U{211d} }%
%BeginExpansion
\mathbb{R}
%EndExpansion
,\\
\sigma(t,s,x,u) & :\left[  0,T\right]  ^{2}\times%
%TCIMACRO{\U{211d} }%
%BeginExpansion
\mathbb{R}
%EndExpansion
\times\mathbb{U}\times\Omega\rightarrow%
%TCIMACRO{\U{211d} }%
%BeginExpansion
\mathbb{R}
%EndExpansion
,\\
g(t,s,x,y,z,k(\cdot),u) & :\left[  0,T\right]  ^{2}\times%
%TCIMACRO{\U{211d} }%
%BeginExpansion
\mathbb{R}
%EndExpansion
^{3}\times L^{2}(\nu)\times\mathbb{U}\times\Omega\rightarrow%
%TCIMACRO{\U{211d} }%
%BeginExpansion
\mathbb{R}
%EndExpansion
,\\
\theta(t,s,x,u,e) & :\left[  0,T\right]  ^{2}\times%
%TCIMACRO{\U{211d} }%
%BeginExpansion
\mathbb{R}
%EndExpansion
\times\mathbb{U}\times%
%TCIMACRO{\U{211d} }%
%BeginExpansion
\mathbb{R}
%EndExpansion
_{0}\times\Omega\rightarrow%
%TCIMACRO{\U{211d} }%
%BeginExpansion
\mathbb{R}
%EndExpansion
,
\end{array}
\]
are assumed to be $C^{1}$ with respect to their first variables, and for all
$t,x,y,z,k,u,e$ the processes $s\mapsto b(t,s,x,u),s\mapsto\sigma
(t,s,x,u),s\mapsto g(t,s,x,y,z,k(\cdot),u),s\mapsto\theta(t,s,x,u,e)$ are
$\mathcal{F}_{s}$-measurable for all $s\leq t$. We assume that $t\mapsto
Z\left(  t,s\right)  $ and $t\mapsto K\left(  t,s,\cdot\right)  $ are $C^{1}$
for all $s,e,\omega$ and that%
\begin{equation}
\mathbb{E}\left[  \int_{0}^{T}\int_{0}^{T}\left(  \frac{\partial Z}{\partial
t}\left(  t,s\right)  \right)  ^{2}dsdt+\int_{0}^{T}\int_{0}^{T}%
\int_{\mathbb{R}_{0}}\left(  \frac{\partial K}{\partial t}\left(
t,s,e\right)  \right)  ^{2}\nu\left(  de\right)  dsdt\right]  <\infty
.\label{eq2.3a}%
\end{equation}
It is known that \eqref{eq2.3a} holds for some linear systems. See
\cite{HO}.\newline Let $\mathbb{U}$ be a given open convex subset of
$\mathbb{R}$ and let $\mathcal{U}=\mathcal{U}_{\mathbb{G}}$ be a given family
of \emph{admissible} controls, required to be $\mathbb{G}-$predictable, where,
as before, $\mathbb{G=}\{\mathcal{G}_{t}\}_{t\geq0}$ is a given subfiltration
of $\mathbb{F=}\{\mathcal{F}_{t}\}_{t\geq0}$, in the sense that $\mathcal{G}%
_{t}\subseteq\mathcal{F}_{t}$ for all $t.$
%and we denote
%by $\mathcal{L}$ the set of all stochastic processes with parameter space
%$[0,T]$.\newline
We associate to the system $\left(  \ref{a1}\right)  -\left(  \ref{a2}\right)
$ the following \emph{performance functional}:
\begin{equation}
J(u)=\mathbb{E}\left[  \int_{0}^{T}f(s,X(s),Y(s),u(s))ds+\varphi
(X(T))+\psi(Y(0))\right]  ,\label{J}%
\end{equation}
for given functions
\begin{align*}
f &  :\left[  0,T\right]  \times\mathbb{R}^{2}\times\mathbb{U}\times
\Omega\rightarrow\mathbb{R},\\
\varphi &  :\mathbb{R}\rightarrow\mathbb{R},\\
\psi &  :\mathbb{R}\rightarrow\mathbb{R}.
\end{align*}
The functions $\varphi,\psi$ are assumed to be $C^{1}$, while $f(s,x,y,u)$ is
assumed to be $\mathbb{F}$-adapted with respect to $s$ and $C^{1}$ with
respect to $x,y,u$ for each s. We remark here that our performance functional
is not of Volterra type. Our optimisation control problem is to find $u^{\ast
}\in\mathcal{U}_{\mathbb{G}}$ such that
\begin{equation}
\underset{u\in\mathcal{U}}{\text{sup}}\text{ }J(u)=J(u^{\ast}).\label{j}%
\end{equation}
Let $\mathcal{L}$ be the set of all $\mathbb{F}$-adapted stochastic processes,
and let $\mathcal{R}$ denote the set of all functions $k:$ $\mathbb{R}%
_{0}\rightarrow\mathbb{R}.$ Define the \emph{Hamiltonian functional}:%
\begin{align}
\mathcal{H}(t,x,y,z,k\left(  \cdot\right)  ,v,p,p(\cdot),q,\lambda
,\lambda(\cdot),r\left(  \cdot\right)  ) &  :=H_{0}(t,x,y,z,k\left(
\cdot\right)  ,v,p,q,\lambda,r\left(  \cdot\right)  )\label{eq3.3}\\
&  +H_{1}(t,x,y,z,k\left(  \cdot\right)  ,v,p(\cdot),\lambda(\cdot)),\nonumber
\end{align}
where%

\[
H_{0}:[0,T]\times\mathbb{R}^{3}\mathbb{\times\mathcal{R\times}U\times R}%
^{3}\times\mathbb{\mathcal{R}}\rightarrow\mathbb{R}%
\]
is defined by
\begin{equation}%
\begin{array}
[c]{l}%
H_{0}(t,x,y,z,k\left(  \cdot\right)  ,v,p,q,\lambda,r\left(  \cdot\right)  )\\
:=f(t,x,y,v)+b(t,t,x,v)p+\sigma(t,t,x,v)q\\
+\int_{%
%TCIMACRO{\U{211d} }%
%BeginExpansion
\mathbb{R}
%EndExpansion
_{0}}\theta(t,t,x,v)r(t,e)\nu(de)+g(t,t,x,y,z,k\left(  \cdot\right)
,v)\lambda
\end{array}
\label{eq2.1}%
\end{equation}
and%
\[
H_{1}:[0,T]\times\mathbb{R}^{3}\mathbb{\times}\mathcal{R}\times\mathbb{U}%
\times\mathcal{L}\times\mathcal{L}\rightarrow\mathbb{R}%
\]
is defined by%
\begin{align}
&  H_{1}(t,x,y,z,k\left(  \cdot\right)  ,v,p(\cdot),\lambda(\cdot
))\label{eq2.8}\\
&  :=\int_{t}^{T}\dfrac{\partial b}{\partial s}(s,t,x,v)p(s)ds+\int_{t}%
^{T}\dfrac{\partial\sigma}{\partial s}(s,t,x,v)\mathbb{E}[\left.
D_{t}p(s)\right\vert \mathcal{F}_{t}]ds\nonumber\\
&  +\int_{t}^{T}\int_{\mathbb{R}_{0}}\dfrac{\partial\theta}{\partial
s}(s,t,x,v)\mathbb{E}[\left.  D_{t,e}p(s)\right\vert \mathcal{F}_{t}%
]\nu(de)ds+\int_{0}^{t}\dfrac{\partial g}{\partial s}(s,t,x,y,z,k\left(
\cdot\right)  ,v)\lambda(s)ds\nonumber\\
&  +\int_{0}^{t}\dfrac{\partial g}{\partial z}(s,t,x,y,z,k\left(
\cdot\right)  ,v)\frac{\partial Z}{\partial s}(s,t)\lambda(s)ds\nonumber\\
&  +\int_{0}^{t}\left\langle \nabla_{k}g(s,t,x,y,z,k\left(  \cdot\right)
,v),\frac{\partial K}{\partial s}(s,t,\cdot)\right\rangle \lambda
(s)ds.\nonumber
\end{align}
Here, and in the following, $D_{t}$ and $D_{t,e}$ denote the
\emph{(generalised) Hida-Malliavin derivatives} at $t$ and at $(t,e)$ with
respect to $B$ and $\tilde{N}$, respectively, and $\nabla_{k}$ denotes the
\emph{Fr\'{e}chet derivative} with respect to $k$. We refer to the Appendix
for more details.\newline\newline The associated \emph{forward-backward
system} for the adjoint processes $\lambda\left(  t\right)  $, \newline%
$\left(  p\left(  t\right)  ,q\left(  t\right)  ,r(t,\cdot)\right)  $ is
\begin{equation}%
\begin{cases}
d\lambda(t):=\frac{\partial\mathcal{H}}{\partial y}(t)dt+\frac{\partial
\mathcal{H}}{\partial z}(t)dB(t)+\int_{%
%TCIMACRO{\U{211d} }%
%BeginExpansion
\mathbb{R}
%EndExpansion
_{0}}\frac{d\nabla_{k}\mathcal{H}}{d\nu}(t)\tilde{N}\left(  dt,de\right)
,0\leq t\leq T,\\
\lambda(0):=\psi^{\prime}\left(  Y\left(  0\right)  \right)  ,
\end{cases}
\label{eq3.2}%
\end{equation}
and%

\begin{equation}%
\begin{cases}
dp(t):=-\frac{\partial\mathcal{H}}{\partial x}(t)dt+q(t)dB(t)+\int_{%
%TCIMACRO{\U{211d} }%
%BeginExpansion
\mathbb{R}
%EndExpansion
_{0}}r(t,e)\tilde{N}\left(  dt,de\right)  ,0\leq t\leq T,\\
p(T):=\varphi^{\prime}(X(T))+\lambda\left(  T\right)  \eta^{\prime}\left(
X\left(  T\right)  \right)  ,
\end{cases}
\label{eq3.4}%
\end{equation}
where we have used the simplified notation
\begin{equation}
\frac{\partial\mathcal{H}}{\partial x}(t)=\left[  \frac{\partial\mathcal{H}%
}{\partial x}(t,x,Y\left(  t\right)  ,Z\left(  t,.\right)  ,K(t,\cdot
),u\left(  t\right)  ,p(t),q(t),\lambda\left(  t\right)  ,r(t,\cdot))\right]
_{x=X\left(  t\right)  },\label{eq3.5}%
\end{equation}
and similarly for $\frac{\partial\mathcal{H}}{\partial y}(t)$, $\frac
{\partial\mathcal{H}}{\partial z}(t)...$\newline As in \cite{os} we assume
that $H$ is Fr\'{e}chet differentiable $(C^{1})$ in the variables $x,y,z,k,u$
and that the Fr\'{e}chet derivative $\nabla_{k}H$ of $H$ with respect to
$k\in\mathcal{R}$ as a random measure is absolutely continuous with respect to
$\nu$, with Radon-Nikodym derivative $\displaystyle\frac{d\nabla_{k}H}{d\nu}$.
Thus, if $\langle\nabla_{k}H,h\rangle$ denotes the action of the linear
operator $\nabla_{k}H$ on the function $h\in\mathcal{R},$ we have
\begin{equation}
\langle\nabla_{k}H,h\rangle=\int_{\mathbb{R}_{0}}h(e)d\nabla_{k}%
H(e)=\int_{\mathbb{R}_{0}}h(e)\frac{d\nabla_{k}H(e)}{d\nu(e)}d\nu
(e).\label{3.10}%
\end{equation}
The question of existence and uniqueness of the forward-backward system above
will not be studied here. It is a subject of future research. See, however our
partial result in Section 3.

\subsection{A sufficient maximum principle}

In this subsection, we prove that under some conditions such as the concavity,
a given control $\hat{u}$ which satisfies a maximum condition of the
Hamiltonian, is an optimal control for the problem $\left(  \ref{j}\right)
$.\newline From $\left(  \ref{a1}\right)  -\left(  \ref{a2}\right)  $ we can
get the differential forms:
\begin{equation}%
\begin{array}
[c]{l}%
dX(t)=\xi^{\prime}(t)dt+b\left(  t,t,X(t),u(t)\right)  dt+\left(
%TCIMACRO{\dint _{0}^{t}}%
%BeginExpansion
{\displaystyle\int_{0}^{t}}
%EndExpansion
\dfrac{\partial b}{\partial t}\left(  t,s,X(s),u(s)\right)  ds\right)  dt\\
+\sigma\left(  t,t,X(t),u(t)\right)  dB(t)+\left(
%TCIMACRO{\dint _{0}^{t}}%
%BeginExpansion
{\displaystyle\int_{0}^{t}}
%EndExpansion
\dfrac{\partial\sigma}{\partial t}\left(  t,s,X(s),u(s)\right)  dB(s)\right)
dt\\
+%
%TCIMACRO{\dint _{\mathbb{R}_{0}}}%
%BeginExpansion
{\displaystyle\int_{\mathbb{R}_{0}}}
%EndExpansion
\theta(t,t,X(t),u(t),e)\tilde{N}(dt,de)+\left(
%TCIMACRO{\dint _{0}^{t}}%
%BeginExpansion
{\displaystyle\int_{0}^{t}}
%EndExpansion%
%TCIMACRO{\dint _{\mathbb{R}_{0}}}%
%BeginExpansion
{\displaystyle\int_{\mathbb{R}_{0}}}
%EndExpansion
\dfrac{\partial\theta}{\partial t}(t,s,X(s),u(s),e)\tilde{N}(ds,de)\right)
dt,
\end{array}
\label{eq3.6}%
\end{equation}
and
\begin{align}
dY(t) &  =-g(t,t,X(t),Y(t),Z(t,t),K(t,t,\cdot),u(t))dt\nonumber\\
&  +\left(  \int_{t}^{T}\frac{\partial g}{\partial t}%
(t,s,X(s),Y(s),Z(t,s),K(t,s,\cdot),u(s))ds\right)  dt\nonumber\\
&  +\int_{t}^{T}\frac{\partial g}{\partial z}(t,s,X(s),Y(s),Z(t,s),K(t,s,\cdot
),u(s))\frac{\partial Z}{\partial t}\left(  t,s\right)  dt\nonumber\\
&  +\int_{t}^{T}\left\langle \nabla_{k}g(t,s,X(s),Y(s),Z(t,s),K(t,s,\cdot
),u(s)),\frac{\partial K}{\partial t}\left(  t,s,\cdot\right)  \right\rangle
dt\nonumber\\
&  +Z(t,t)dB(t)+\int_{\mathbb{R}_{0}}K(t,t,e)\tilde{N}(dt,de)\nonumber\\
&  -\left(  \int_{t}^{T}\frac{\partial Z}{\partial t}(t,s)dB(s)\right)
dt-\left(  \int_{t}^{T}\int_{\mathbb{R}_{0}}\frac{\partial K}{\partial
t}(t,s,e)\tilde{N}(ds,de)\right)  dt,\nonumber\\
Y(T) &  =\eta(X(T)).\label{eq2.13a}%
\end{align}
We now state and prove a sufficient maximum principle:

\begin{theorem}
Let $\hat{u}\in\mathcal{U}_{\mathbb{G}},$ with corresponding solutions
$\hat{X}(t),\newline(\hat{Y}(t),\hat{Z}(t,s),\hat{K}\left(  t,s,\cdot\right)
),\hat{\lambda}(t),\left(  \hat{p}(t),\hat{q}(t),\hat{r}\left(  t,\cdot
\right)  \right)  $ of equations
\eqref {eq3.6},\eqref {eq2.13a},\eqref {eq3.2} and \eqref{eq3.4},
respectively. Assume the following:$\cdot$

\begin{itemize}
\item (Concavity conditions) The functions%
\[
x\mapsto\eta(x),\text{ }x\mapsto\varphi\left(  x\right)  ,\text{ }x\mapsto
\psi\left(  x\right)  \text{ }%
\]
and
\[
x,y,z,k\left(  \cdot\right)  ,u\mapsto\mathcal{H}(t,x,y,z,k\left(
\cdot\right)  ,u,p,q,\lambda,r),
\]
are concave for all $t,p,q,\lambda,r$.\newline

\item (The maximum condition)
\begin{align}
&  \underset{v\in\mathcal{U}}{\sup}\text{ }\mathbb{E}\left[  \mathcal{H}%
(t,\hat{X}\left(  t\right)  ,\hat{Y}\left(  t\right)  ,\hat{Z}\left(
t\right)  ,\hat{k}\left(  t,\cdot\right)  ,v,\hat{\lambda}\left(  t\right)
,\hat{p}(t),\hat{q}(t),\hat{r}(t,\cdot))\mid{\mathcal{G}}_{t}\right]
\nonumber\\
&  =\mathbb{E}\left[  \mathcal{H}(t,\hat{X}(t),\hat{Y}\left(  t\right)
,\hat{Z}\left(  t\right)  ,\hat{k}\left(  t,\cdot\right)  ,\hat{u}%
(t),\hat{\lambda}\left(  t\right)  ,\hat{p}(t),\hat{q}(t),\hat{r}%
(t,\cdot))\mid{\mathcal{G}}_{t}\right]  ,\text{ }\forall t\geq0. \label{eq3.9}%
\end{align}
%and that%
%\begin{equation}
%\text{(Growth conditions)} \label{1.8}%
%\end{equation}

\end{itemize}

Then, $\hat{u}$ is an optimal $\mathbb{G}-$adapted control.
\end{theorem}

\noindent{Proof.} \quad By considering a suitable increasing family of
stopping times converging to $T$, we may assume that all the local martingales
appearing in the proof below are martingales. In particular, the expectations
of the $dB$- and $\tilde{N}(dt,de)$-integrals are all $0$.\newline Choose an
arbitrary $u\in$ $\mathcal{U}_{\mathbb{G}}$ and consider
\begin{equation}
J(u)-J(\hat{u})=I_{1}+I_{2}+I_{3},\nonumber
\end{equation}
where
\begin{align}
I_{1} &  =\mathbb{E}\left[
%TCIMACRO{\dint _{0}^{T}}%
%BeginExpansion
{\displaystyle\int_{0}^{T}}
%EndExpansion
\left\{  f\left(  t\right)  -\hat{f}\left(  t\right)  \right\}  dt\text{
}\right]  ,\quad I_{2}=\mathbb{E}\left[  \varphi\left(  X(T\right)
)-\varphi(\hat{X}\left(  T\right)  )\right]  ,\text{ }\nonumber\\
I_{3} &  =\mathbb{E}\left[  \psi\left(  Y\left(  0\right)  \right)  -\psi
(\hat{Y}\left(  0\right)  )\right]  ,\label{eq2.8}%
\end{align}
where $f(t)=$ $f\left(  t,X(t),Y\left(  t\right)  ,u(t)\right)  ,$ $\hat
{f}\left(  t\right)  =f(t,\hat{X}(t),\hat{Y}\left(  t\right)  ,\hat{u}(t)).$
Using a simplified notation
\[%
\begin{array}
[c]{l}%
b(t,t)=b\left(  t,t,X(t),u(t)\right)  ,\hat{b}\left(  t,t\right)
=b(t,t,\hat{X}(t),\hat{u}(t)),b(t,s)=b\left(  t,s,X(s),u(s)\right)  \\
\theta(t,t,e)=\theta(t,t,X(s),u(s),e),\theta(t,s,e)=\theta
(t,s,X(s),u(s),e)\text{ etc., we get}%
\end{array}
\text{ }%
\]%
\begin{align}
I_{1} &  =\mathbb{E}\left[
%TCIMACRO{\dint _{0}^{T}}%
%BeginExpansion
{\displaystyle\int_{0}^{T}}
%EndExpansion
\left\{  H_{0}(t)-\hat{H}_{0}(t)-\hat{p}(t)\left(  b(t,t)-\hat{b}(t,t)\right)
-\hat{q}(t)\left(  \sigma(t,t)-\hat{\sigma}(t,t)\right)  \right.  \right.
\nonumber\\
&  \left.  \left.  -\hat{\lambda}\left(  t\right)  \left(  g\left(
t,t\right)  -\hat{g}\left(  t,t\right)  \right)  -%
%TCIMACRO{\dint _{\mathbb{R}_{0}}}%
%BeginExpansion
{\displaystyle\int_{\mathbb{R}_{0}}}
%EndExpansion
\hat{r}\left(  t,e\right)  \left(  \theta\left(  t,t,e\right)  -\hat{\theta
}\left(  t,t,e\right)  \right)  \nu\left(  de\right)  \right\}  dt\right]
.\label{eq2.9}%
\end{align}
Using concavity and the It\^{o} formula, we obtain
\begin{align}
I_{2} &  \leq\mathbb{E}\left[  \varphi^{\prime}(\hat{X}(T))\left(
X(T)-\hat{X}(T)\right)  \right]  \nonumber\\
&  =\mathbb{E}\left[  \hat{p}(T)\left(  X(T)-\hat{X}(T)\right)  \right]
-\mathbb{E}\left[  \hat{\lambda}\left(  T\right)  \eta^{\prime}(\hat{X}\left(
T\right)  )\left(  X(T)-\hat{X}(T)\right)  \right]  \nonumber\\
&  =\mathbb{E}\left[  \int_{0}^{T}\hat{p}(t)\left(  dX(t)-d\hat{X}(t)\right)
+\int_{0}^{T}\left(  X(t)-\hat{X}(t)\right)  d\hat{p}(t)\right.  \nonumber\\
&  +\int_{0}^{T}\hat{q}(t)(\sigma(t,t)-\hat{\sigma}(t,t))dt+\int_{0}^{T}%
\int_{\mathbb{R}_{0}}\hat{r}(t,e)(\theta(t,t,e)-\hat{\theta}(t,t,e))\nu
(de)dt]\nonumber\\
&  -\mathbb{E}\left[  \hat{\lambda}\left(  T\right)  \eta^{\prime}(\hat
{X}\left(  T\right)  )\left(  X(T)-\hat{X}(T)\right)  \right]  \nonumber\\
&  =\mathbb{E}[\int_{0}^{T}\{\hat{p}(t)\left(  b(t,t)-\hat{b}(t,t)+\int
_{0}^{t}(\frac{\partial b}{\partial t}(t,s)-\frac{\partial\hat{b}}{\partial
t}(t,s))ds\right.  \nonumber\\
&  +\int_{0}^{t}\left(  \frac{\partial\sigma}{\partial t}(t,s)-\frac
{\partial\hat{\sigma}}{\partial t}(t,s)\right)  dB(s)\nonumber\\
&  \left.  +\int_{0}^{t}\int_{\mathbb{R}_{0}}\left(  \frac{\partial\theta
}{\partial t}(t,s,e)-\frac{\partial\hat{\theta}}{\partial t}(t,s,e)\right)
\tilde{N}(ds,de)\right)  -\frac{\partial\hat{\mathcal{H}}}{\partial
x}(t)(X(t)-\hat{X}(t))\nonumber\\
&  +\hat{q}(t)[\sigma(t,t)-\hat{\sigma}(t,t)]\}dt]+\int_{0}^{T}\int
_{\mathbb{R}_{0}}\hat{r}(t,e)(\theta(t,t,e)-\hat{\theta}(t,t,e))\nu
(de)dt]\nonumber\\
&  -\mathbb{E}[\hat{\lambda}(T)\eta^{\prime}(\hat{X}(T))(X(T)-\hat
{X}(T))].\label{eq2.10}%
\end{align}
By the Fubini theorem, we get
\begin{align}
\int_{0}^{T}\left(  \int_{0}^{t}\frac{\partial b}{\partial t}(t,s)ds\right)
\hat{p}(t)dt &  =\int_{0}^{T}\left(  \int_{s}^{T}\frac{\partial b}{\partial
t}(t,s)\hat{p}(t)dt\right)  ds\nonumber\\
&  =\int_{0}^{T}\left(  \int_{t}^{T}\frac{\partial b}{\partial s}(s,t)\hat
{p}(s)ds\right)  dt,\label{eq2.11}%
\end{align}
and by the generalised duality theorems for the Malliavin derivatives
\cite{AO}, we have%
\begin{align}
\mathbb{E}\left[  \int_{0}^{T}\left(  \int_{0}^{t}\frac{\partial\sigma
}{\partial t}(t,s)dB(s)\right)  \hat{p}(t)dt\right]   &  =\int_{0}%
^{T}\mathbb{E}\left[  \int_{0}^{t}\frac{\partial\sigma}{\partial
t}(t,s)dB(s)\hat{p}(t)\right]  dt\nonumber\\
&  =\int_{0}^{T}\mathbb{E}\left[  \int_{0}^{t}\frac{\partial\sigma}{\partial
t}(t,s)\mathbb{E}[D_{s}\hat{p}(t)\mid\mathcal{F}_{s}]ds\right]  dt\nonumber\\
&  =\int_{0}^{T}\mathbb{E}\left[  \int_{s}^{T}\frac{\partial\sigma}{\partial
t}(t,s)\mathbb{E}[D_{s}\hat{p}(t)\mid\mathcal{F}_{s}]dt\right]  ds\nonumber\\
&  =\mathbb{E}\left[  \int_{0}^{T}\int_{t}^{T}\frac{\partial\sigma}{\partial
s}(s,t)\mathbb{E}[D_{t}\hat{p}(s)\mid\mathcal{F}_{t}]dsdt\right]
,\label{eq2.12}%
\end{align}
and%
\begin{align}
&  \mathbb{E}\left[  \int_{0}^{T}\left(  \int_{0}^{t}\int_{%
%TCIMACRO{\U{211d} }%
%BeginExpansion
\mathbb{R}
%EndExpansion
_{0}}\left(  \frac{\partial\theta}{\partial t}(t,s,e)\right)  \tilde
{N}(ds,de)p\left(  t\right)  \right)  dt\right]  \nonumber\\
&  =\int_{0}^{T}\mathbb{E}\left[  \int_{0}^{t}\int_{%
%TCIMACRO{\U{211d} }%
%BeginExpansion
\mathbb{R}
%EndExpansion
_{0}}\left(  \frac{\partial\theta}{\partial t}(t,s,e)\right)  \tilde
{N}(ds,de)p\left(  t\right)  \right]  dt\nonumber\\
&  =\int_{0}^{T}\mathbb{E}\left[  \int_{s}^{T}\int_{%
%TCIMACRO{\U{211d} }%
%BeginExpansion
\mathbb{R}
%EndExpansion
_{0}}\frac{\partial\theta}{\partial t}(t,s,e)\mathbb{E}\left[  \left.
D_{s,e}p\left(  t\right)  \right\vert \mathcal{F}_{s}\right]  \nu\left(
de\right)  dt\right]  ds\nonumber\\
&  =\int_{0}^{T}\mathbb{E}\left[  \int_{t}^{T}\int_{%
%TCIMACRO{\U{211d} }%
%BeginExpansion
\mathbb{R}
%EndExpansion
_{0}}\frac{\partial\theta}{\partial s}(s,t,e)\mathbb{E}\left[  \left.
D_{t,e}p\left(  s\right)  \right\vert \mathcal{F}_{t}\right]  \nu\left(
de\right)  ds\right]  dt.\label{eq2.17}%
\end{align}
Substituting \eqref{eq2.12}, \eqref{eq2.17} and \eqref{eq2.11} into
\eqref{eq2.10}, we get%
\begin{align}
I_{2} &  \left.  \leq\mathbb{E}\left[
%TCIMACRO{\dint _{0}^{T}}%
%BeginExpansion
{\displaystyle\int_{0}^{T}}
%EndExpansion
\left\{  \hat{p}(t)\left(  b(t,t)-\hat{b}\left(  t,t\right)  \right)
+\int_{t}^{T}\hat{p}(s)\left(  \frac{\partial b}{\partial s}\left(
s,t\right)  -\frac{\partial\hat{b}}{\partial s}\left(  s,t\right)  \right)
ds\right.  \right.  \right.  \nonumber\\
&  +\int_{t}^{T}\left(  \frac{\partial\sigma}{\partial s}(s,t)-\frac
{\partial\hat{\sigma}}{\partial s}(s,t)\right)  \mathbb{E}[D_{t}\hat{p}%
(s)\mid\mathcal{F}_{t}]ds\nonumber\\
&  +\int_{t}^{T}\int_{%
%TCIMACRO{\U{211d} }%
%BeginExpansion
\mathbb{R}
%EndExpansion
_{0}}\frac{\partial\theta}{\partial s}(s,t,e)\mathbb{E}\left[  \left.
D_{t,e}p\left(  s\right)  \right\vert \mathcal{F}_{t}\right]  \nu\left(
de\right)  ds\nonumber\\
&  +\int\limits_{0}^{T}\int\limits_{%
%TCIMACRO{\U{211d} }%
%BeginExpansion
\mathbb{R}
%EndExpansion
_{0}}\left(  \hat{r}\left(  t,e\right)  \left(  \theta\left(  t,t,e\right)
-\hat{\theta}\left(  t,t,e\right)  \right)  \nu\left(  de\right)  \right)
dt\nonumber\\
&  \left.  \left.  -\frac{\partial\hat{\mathcal{H}}}{\partial x}(t)\left(
X(t)-\hat{X}(t)\right)  +\hat{q}(t)\left(  \sigma\left(  t,t\right)
-\hat{\sigma}\left(  t,t\right)  \right)  \right\}  dt\right]  \nonumber\\
&  -\mathbb{E}\left[  \hat{\lambda}\left(  T\right)  \eta^{\prime}(\hat
{X}\left(  T\right)  )\left(  X(T)-\hat{X}(T)\right)  \right]  .\label{eq2.13}%
\end{align}
By the concavity of $\psi$ and $\eta$, we obtain%
\begin{align}
I_{3} &  =\mathbb{E}\left[  \psi\left(  Y\left(  0\right)  \right)
-\psi\left(  \hat{Y}\left(  0\right)  \right)  \right]  \nonumber\\
&  \leq\mathbb{E}\left[  \psi^{\prime}\left(  \hat{Y}\left(  0\right)
\right)  \left(  Y\left(  0\right)  -\hat{Y}\left(  0\right)  \right)
\right]  \nonumber\\
&  =\mathbb{E}\left[  \hat{\lambda}\left(  0\right)  \left(  Y\left(
0\right)  -\hat{Y}\left(  0\right)  \right)  \right]  \nonumber
\end{align}

\begin{align}
&  =\mathbb{E}\left[  \hat{\lambda}\left(  T\right)  \left(  Y\left(
T\right)  -\hat{Y}\left(  T\right)  \right)  \right]  -\mathbb{E}\left[
\int_{0}^{T}\left(  Y\left(  t\right)  -\hat{Y}\left(  t\right)  \right)
d\hat{\lambda}\left(  t\right)  \right.  \nonumber\\
&  \left.  +\int_{0}^{T}\hat{\lambda}\left(  t\right)  \left(  dY\left(
t\right)  -d\hat{Y}\left(  t\right)  \right)  +\int_{0}^{T}\frac{\partial
\hat{\mathcal{H}}}{\partial z}\left(  t\right)  \left(  Z\left(  t,s\right)
-\hat{Z}\left(  t,s\right)  \right)  dt\right.  \nonumber\\
&  \left.  +\int_{0}^{T}\int_{\mathbb{R}_{0}}\frac{d\nabla_{k}\mathcal{H}%
}{d\nu}(t)(K(t,s,e)-\hat{K}(t,s,e))\nu(de)dt\right]  \nonumber\\
&  \leq\mathbb{E}\left[  \hat{\lambda}\left(  T\right)  \eta^{\prime}\left(
X\left(  T\right)  \right)  \left(  X(T)-\hat{X}(T)\right)  \right]
\nonumber\\
&  -\mathbb{E}\left[  \int_{0}^{T}\frac{\partial\hat{\mathcal{H}}}{\partial
y}\left(  t\right)  \left(  Y\left(  t\right)  -\hat{Y}\left(  t\right)
\right)  dt-\int_{0}^{T}\hat{\lambda}\left(  t\right)  \left(  g\left(
t,t\right)  -\hat{g}\left(  t,t\right)  \right)  dt\right.  \nonumber\\
&  \left.  +\int_{0}^{T}\left(  \hat{\lambda}\left(  t\right)  \int_{t}%
^{T}\left(  \frac{\partial g}{\partial t}\left(  t,s\right)  -\frac
{\partial\hat{g}}{\partial t}\left(  t,s\right)  \right)  ds\right)
dt\right.  \nonumber\\
&  +\int_{0}^{T}\hat{\lambda}\left(  t\right)  \left[  \int_{t}^{T}\left(
\frac{\partial g}{\partial z}\left(  t,s\right)  \frac{\partial Z}{\partial
t}\left(  t,s\right)  -\frac{\partial\hat{g}}{\partial z}\left(  t,s\right)
\frac{\partial\hat{Z}}{\partial t}\left(  t,s\right)  \right)  ds\right]
dt\nonumber\\
&  +\int_{0}^{T}\hat{\lambda}\left(  t\right)  \left[  \int_{t}^{T}\left(
\left\langle \nabla_{k}g\left(  t,s\right)  ,\frac{\partial K}{\partial
t}\left(  t,s,\cdot\right)  \right\rangle -\left\langle \nabla_{k}\hat
{g}\left(  t,s\right)  ,\frac{\partial\hat{K}}{\partial t}\left(
t,s,\cdot\right)  \right\rangle \right)  ds\right]  dt\nonumber\\
&  \left.  +\int_{0}^{T}\left(  \hat{\lambda}\left(  t\right)  \int_{t}%
^{T}\left(  \frac{\partial Z}{\partial t}\left(  t,s\right)  -\frac
{\partial\hat{Z}}{\partial t}\left(  t,s\right)  \right)  dB(s)\right)
dt\right.  \nonumber\\
&  \left.  +\int_{0}^{T}\left(  \hat{\lambda}\left(  t\right)  \int_{t}%
^{T}\int_{\mathbb{R}_{0}}\left(  \frac{\partial K}{\partial t}\left(
t,s,\cdot\right)  -\frac{\partial\hat{K}}{\partial t}\left(  t,s,\cdot\right)
\right)  \tilde{N}(ds,de)\right)  dt\right.  \nonumber\\
&  \left.  +\int_{0}^{T}\frac{\partial\hat{\mathcal{H}}}{\partial z}\left(
t\right)  \left(  Z\left(  t,s\right)  -\hat{Z}\left(  t,s\right)  \right)
dt\right.  \nonumber\\
&  \left.  +\int_{0}^{T}\int_{\mathbb{R}_{0}}\frac{d\nabla_{k}\hat
{\mathcal{H}}}{d\nu}(t)\left(  K(t,s,e)-\hat{K}(t,s,e)\right)  \nu
(de)dt\right]  .\label{eq2.14}%
\end{align}
By the Fubini Theorem, we get
\begin{align}
\int_{0}^{T}\left(  \int_{t}^{T}\frac{\partial g}{\partial t}(t,s)ds\right)
\hat{\lambda}(t)dt &  =\int_{0}^{T}\left(  \int_{0}^{s}\frac{\partial
g}{\partial t}(t,s)\hat{\lambda}(t)dt\right)  ds\nonumber\\
&  =\int_{0}^{T}\left(  \int_{0}^{t}\frac{\partial g}{\partial s}%
(s,t)\hat{\lambda}(s)ds\right)  dt,\label{eq2.15}%
\end{align}%
\begin{equation}%
\begin{array}
[c]{l}%
\int_{0}^{T}\hat{\lambda}(t)\left[  \int_{t}^{T}\frac{\partial g}{\partial
z}(t,s)\frac{\partial Z}{\partial t}\left(  t,s\right)  ds\right]  dt=\int
_{0}^{T}\left(  \int_{0}^{t}\hat{\lambda}(s)\frac{\partial g}{\partial
z}(s,t)\frac{\partial Z}{\partial s}\left(  s,t\right)  ds\right)  dt,
\end{array}
\label{eq2.16}%
\end{equation}
and%
\begin{equation}%
\begin{array}
[c]{l}%
%TCIMACRO{\dint _{0}^{T}}%
%BeginExpansion
{\displaystyle\int_{0}^{T}}
%EndExpansion
\hat{\lambda}(t)\left[
%TCIMACRO{\dint _{t}^{T}}%
%BeginExpansion
{\displaystyle\int_{t}^{T}}
%EndExpansion
\left\langle \nabla_{k}g\left(  t,s\right)  ,\dfrac{\partial K}{\partial
t}\left(  t,s,\cdot\right)  \right\rangle ds\right]  dt\\
=%
%TCIMACRO{\dint _{0}^{T}}%
%BeginExpansion
{\displaystyle\int_{0}^{T}}
%EndExpansion
\left(
%TCIMACRO{\dint _{0}^{t}}%
%BeginExpansion
{\displaystyle\int_{0}^{t}}
%EndExpansion
\hat{\lambda}(s)\left\langle \nabla_{k}g\left(  s,t\right)  ,\dfrac{\partial
K}{\partial s}\left(  s,t,\cdot\right)  \right\rangle ds\right)  dt.
\end{array}
\label{eq.2.19}%
\end{equation}
Substituting \eqref{eq2.15}-\eqref{eq.2.19} into \eqref{eq2.14}, we get%
\begin{align}
I_{3} &  \leq\mathbb{E}\left[  \hat{\lambda}\left(  T\right)  \eta^{\prime
}\left(  X\left(  T\right)  \right)  \left(  X(T)-\hat{X}(T)\right)  \right]
\nonumber\\
&  -\mathbb{E}\left[  \int_{0}^{T}\frac{\partial\hat{\mathcal{H}}}{\partial
y}\left(  t\right)  \left(  Y\left(  t\right)  -\hat{Y}\left(  t\right)
\right)  dt-\int_{0}^{T}\hat{\lambda}\left(  t\right)  \left(  g\left(
t,t\right)  -\hat{g}\left(  t,t\right)  \right)  dt\right.  \nonumber\\
&  +\int_{0}^{T}\int_{0}^{t}\left(  \frac{\partial g}{\partial s}\left(
s,t\right)  -\frac{\partial\hat{g}}{\partial s}\left(  s,t\right)  \right)
\hat{\lambda}(s)dsdt\nonumber\\
&  +\int_{0}^{T}\left(  \int_{0}^{t}\hat{\lambda}(s)\left[  \frac{\partial
g}{\partial z}(s,t)\frac{\partial Z}{\partial s}\left(  s,t\right)
-\frac{\partial\hat{g}}{\partial z}(s,t)\frac{\partial\hat{Z}}{\partial
s}\left(  s,t\right)  \right]  ds\right)  dt\nonumber\\
&  +\int_{0}^{T}\left(  \int_{0}^{t}\hat{\lambda}(s)\left[  \left\langle
\nabla_{k}g\left(  s,t\right)  ,\frac{\partial K}{\partial s}\left(
s,t,\cdot\right)  \right\rangle \right.  \right.  \nonumber\\
&  \left.  \left.  \left.  -\left\langle \nabla_{k}\hat{g}\left(  s,t\right)
,\frac{\partial\hat{K}}{\partial s}\left(  s,t,\cdot\right)  \right\rangle
\right]  ds\right)  dt+\int_{0}^{T}\frac{\partial\hat{\mathcal{H}}}{\partial
z}\left(  t\right)  \left(  Z\left(  t,s\right)  -\hat{Z}\left(  t,s\right)
\right)  dt\right]  .\label{eq2.18}%
\end{align}
Adding \eqref{eq2.9}, \eqref{eq2.13} and \eqref{eq2.18}, and noting that%
\[%
\begin{array}
[c]{ll}%
H_{1}\left(  t\right)  -\hat{H}_{1}\left(  t\right)   & =%
%TCIMACRO{\dint _{t}^{T}}%
%BeginExpansion
{\displaystyle\int_{t}^{T}}
%EndExpansion
\left\{  \dfrac{\partial b}{\partial s}\left(  s,t\right)  -\dfrac
{\partial\hat{b}}{\partial s}\left(  s,t\right)  \right\}  \hat{p}(s)ds\\
& +%
%TCIMACRO{\dint _{t}^{T}}%
%BeginExpansion
{\displaystyle\int_{t}^{T}}
%EndExpansion
\left\{  \dfrac{\partial\sigma}{\partial s}(s,t)-\dfrac{\partial\hat{\sigma}%
}{\partial s}(s,t)\right\}  \mathbb{E}[D_{t}\hat{p}(s)\mid\mathcal{F}_{t}]ds\\
& +%
%TCIMACRO{\dint _{0}^{t}}%
%BeginExpansion
{\displaystyle\int_{0}^{t}}
%EndExpansion
\left\{  \dfrac{\partial g}{\partial s}\left(  s,t\right)  -\dfrac
{\partial\hat{g}}{\partial s}\left(  s,t\right)  \right\}  \hat{\lambda
}(s)ds\\
& +%
%TCIMACRO{\dint _{0}^{T}}%
%BeginExpansion
{\displaystyle\int_{0}^{T}}
%EndExpansion
\left(
%TCIMACRO{\dint _{0}^{t}}%
%BeginExpansion
{\displaystyle\int_{0}^{t}}
%EndExpansion
\hat{\lambda}(s)\left[  \dfrac{\partial g}{\partial z}(s,t)\dfrac{\partial
Z}{\partial s}\left(  s,t\right)  -\dfrac{\partial\hat{g}}{\partial
z}(s,t)\dfrac{\partial\hat{Z}}{\partial s}\left(  s,t\right)  \right]
ds\right)  dt\\
& +%
%TCIMACRO{\dint _{0}^{T}}%
%BeginExpansion
{\displaystyle\int_{0}^{T}}
%EndExpansion
\left(
%TCIMACRO{\dint _{0}^{t}}%
%BeginExpansion
{\displaystyle\int_{0}^{t}}
%EndExpansion
\hat{\lambda}(s)\left[  \left\langle \nabla_{k}g\left(  s,t\right)
,\dfrac{\partial K}{\partial s}\left(  s,t,\cdot\right)  \right\rangle
\right.  \right.  \\
& \left.  \left.  -\left\langle \nabla_{k}\hat{g}\left(  s,t\right)
,\frac{\partial\hat{K}}{\partial s}\left(  s,t,\cdot\right)  \right\rangle
\right]  ds\right)  dt\\
& +%
%TCIMACRO{\dint _{t}^{T}}%
%BeginExpansion
{\displaystyle\int_{t}^{T}}
%EndExpansion%
%TCIMACRO{\dint _{\mathbb{R}_{0}}}%
%BeginExpansion
{\displaystyle\int_{\mathbb{R}_{0}}}
%EndExpansion
\left(  \dfrac{\partial\theta}{\partial s}(s,t,e)-\dfrac{\partial\hat{\theta}%
}{\partial s}(s,t,e)\right)  \mathbb{E}[\left.  D_{t,e}p(s)\right\vert
\mathcal{F}_{t}]\nu(de)ds,
\end{array}
\]
we get
\begin{align*}
J(u)-J(\hat{u}) &  =I_{1}+I_{2}+I_{3}\\
&  \leq\mathbb{E}\left[  \int_{0}^{T}\left\{  \mathcal{H}(t)-\hat{\mathcal{H}%
}(t)-\frac{\partial\hat{\mathcal{H}}}{\partial x}(t)\left(  X(t)-\hat
{X}(t)\right)  \right.  \right.  \\
&  \left.  \left.  -\frac{\partial\hat{\mathcal{H}}}{\partial y}\left(
t\right)  \left(  Y\left(  t\right)  -\hat{Y}\left(  t\right)  \right)
-\frac{\partial\hat{\mathcal{H}}}{\partial z}\left(  t\right)  \left(
Z\left(  t,s\right)  -\hat{Z}\left(  t,s\right)  \right)  \right.  \right.  \\
&  \left.  \left.  -\int_{\mathbb{R}_{0}}\frac{d\nabla_{k}\widehat
{\mathcal{H}}}{d\nu}(t)\left(  K(t,s,e)-\hat{K}(t,s,e)\right)  \nu
(de)\right\}  dt\right]  .
\end{align*}
By the concavity of $\mathcal{H}$ and the maximum condition \eqref{eq3.9}, the
proof is complete.

\hfill$\square$ \bigskip

\subsection{A necessary maximum principle}

The concavity condition used in the previous subsection does not always hold
in applications. We prove now if $\hat{u}\in\mathcal{U}_{\mathbb{G}}$ is an
optimal control for the problem $\left(  \ref{j}\right)  $, then we have the
equivalence between being a critical point of $J(u)$ and a critical point of
the conditional Hamiltonian.\newline We start by defining the derivative
processes. For each given $t\in\left[  0,T\right)  $, let $\alpha=\alpha(t)$
be a bounded $\mathcal{G}_{t}-$measurable random variable, let $\epsilon
\in\left(  0, T-t\right]  $ and define%

\begin{equation}
\mu\left(  s\right)  :=\gamma1_{\left[  t,t+\epsilon\right]  }\left(
s\right)  ,\text{ }s\in\left[  0,T\right]  .\label{eq3.1}%
\end{equation}
Assume that
\[
\hat{u}+\epsilon\mu\in\mathcal{U},
\]
for all such $\mu$, and all nonzero $\epsilon$ sufficiently small. Then the
derivative processes are defined by, writing $u$ for $\hat{u}$ for simplicity
from now on,
\[%
\begin{array}
[c]{ll}%
X^{\prime}\left(  t\right)   & :=\dfrac{d}{d\epsilon}X^{u+\epsilon\mu}\left(
t\right)  |_{\epsilon=0},\\
Y^{\prime}\left(  t\right)   & :=\dfrac{d}{d\epsilon}Y^{u+\epsilon\mu}\left(
t\right)  |_{\epsilon=0},\\
Z^{\prime}\left(  t,s\right)   & :=\dfrac{d}{d\epsilon}Z^{u+\epsilon\mu
}\left(  t,s\right)  |_{\epsilon=0},\\
K^{\prime}\left(  t,s,\cdot\right)   & :=\dfrac{d}{d\epsilon}K^{u+\epsilon\mu
}\left(  t,s,\cdot\right)  |_{\epsilon=0}.
\end{array}
\]
We see that%
\begin{align*}
&  X^{\prime}\left(  t\right)  =\int_{0}^{t}\left(  \frac{\partial b}{\partial
x}(t,s)X^{\prime}\left(  s\right)  +\frac{\partial b}{\partial u}%
(t,s)\mu\left(  s\right)  \right)  ds\\
&  +\int_{0}^{t}\left(  \frac{\partial\sigma}{\partial x}(t,s)X^{\prime
}\left(  s\right)  +\frac{\partial\sigma}{\partial u}(t,s)\mu\left(  s\right)
\right)  dB\left(  s\right)  \\
&  +\int_{0}^{t}\int_{\mathbb{R}_{0}}\left(  \frac{\partial\theta}{\partial
x}\left(  t,s,e\right)  X^{\prime}\left(  s\right)  +\frac{\partial\theta
}{\partial u}\left(  t,s,e\right)  \mu\left(  s\right)  \right)  \tilde
{N}(ds,de)
\end{align*}
and
\begin{align*}
Y^{\prime}\left(  t\right)   &  =\eta^{\prime}\left(  X\left(  T\right)
\right)  X^{\prime}\left(  T\right)  +\int_{t}^{T}\left(  \frac{\partial
g}{\partial x}(t,s)X^{\prime}\left(  s\right)  +\frac{\partial g}{\partial
y}(t,s)Y^{\prime}\left(  s\right)  \right.  \\
&  \left.  +\frac{\partial g}{\partial z}(t,s)Z^{\prime}\left(  t,s\right)
+\left\langle \nabla_{k}g\left(  t,s\right)  ,K^{\prime}\left(  t,s,\cdot
\right)  \right\rangle +\frac{\partial g}{\partial u}(t,s)\mu\left(  s\right)
\right)  ds\\
&  -\int_{t}^{T}Z^{\prime}\left(  t,s\right)  dB\left(  s\right)
-\int\nolimits_{t}^{T}\int_{\mathbb{R}_{0}}K^{\prime}\left(  t,s,e\right)
\tilde{N}(ds,de).
\end{align*}
Hence%
\begin{align}
dX^{\prime}\left(  t\right)   &  =\left[  \frac{\partial b}{\partial
x}(t,t)X^{\prime}\left(  t\right)  +\frac{\partial b}{\partial u}%
(t,t)\mu\left(  t\right)  \right.  \nonumber\\
&  +\int_{0}^{t}\left(  \frac{\partial^{2}b}{\partial t\partial x}%
(t,s)X^{\prime}\left(  s\right)  +\frac{\partial^{2}b}{\partial t\partial
u}(t,s)\mu\left(  s\right)  \right)  ds\nonumber\\
&  +\int_{0}^{t}\left(  \frac{\partial^{2}\sigma}{\partial t\partial
x}(t,s)X^{\prime}\left(  s\right)  +\frac{\partial^{2}\sigma}{\partial
t\partial u}(t,s)\mu\left(  s\right)  \right)  dB(s)\nonumber\\
&  \left.  +\int_{0}^{t}\int_{\mathbb{R}_{0}}\left(  \frac{\partial^{2}\theta
}{\partial t\partial x}\left(  t,s,e\right)  X^{\prime}\left(  s\right)
+\frac{\partial^{2}\theta}{\partial t\partial u}\left(  t,s,e\right)
\mu\left(  s\right)  \right)  \tilde{N}(ds,de)\right]  dt\nonumber\\
&  +\left(  \frac{\partial\sigma}{\partial x}(t,t)X^{\prime}\left(  t\right)
+\frac{\partial\sigma}{\partial u}(t,t)\mu\left(  t\right)  \right)
dB(t)\nonumber\\
&  +\int_{%
%TCIMACRO{\U{211d} }%
%BeginExpansion
\mathbb{R}
%EndExpansion
_{0}}\left(  \frac{\partial\theta}{\partial x}\left(  t,t,e\right)  X^{\prime
}\left(  t\right)  +\frac{\partial\theta}{\partial u}\left(  t,t,e\right)
\mu\left(  t\right)  \right)  \tilde{N}(dt,de),\label{eq3.28}%
\end{align}
and%

\begin{equation}%
\begin{array}
[c]{l}%
dY^{\prime}\left(  t\right)  =-\nabla\left(  g(t,t)\right)  \left(  X^{\prime
}\left(  t\right)  ,Y^{\prime}\left(  t\right)  ,Z^{\prime}\left(  t,t\right)
,K^{\prime}\left(  t,t,\cdot\right)  ,\mu\left(  t\right)  \right)  ^{t}dt\\
+\int_{t}^{T}\nabla\left(  \frac{\partial g}{\partial t}(t,s)\right)  \left(
X^{\prime}\left(  s\right)  ,Y^{\prime}\left(  s\right)  ,Z^{\prime}\left(
t,s\right)  ,K^{\prime}\left(  t,s,\cdot\right)  ,\mu\left(  s\right)
\right)  ^{t}dt\\
+\int_{t}^{T}\nabla\left(  \frac{\partial g}{\partial z}(t,s)\right)  \left(
X^{\prime}\left(  s\right)  ,Y^{\prime}\left(  s\right)  ,Z^{\prime}\left(
t,s\right)  ,K^{\prime}\left(  t,s,\cdot\right)  ,\mu\left(  s\right)
\right)  ^{t}\left(  \frac{\partial Z}{\partial t}\left(  t,s\right)  \right)
dt\\
+\int_{t}^{T}\nabla\left(  \nabla_{k}g(t,s)\right)  \left(  X^{\prime}\left(
s\right)  ,Y^{\prime}\left(  s\right)  ,Z^{\prime}\left(  t,s\right)
,K^{\prime}\left(  t,s,\cdot\right)  ,\mu\left(  s\right)  \right)
^{t}\left(  \frac{\partial K}{\partial t}\left(  t,s,\cdot\right)  \right)
dt\\
+Z^{\prime}\left(  t,t\right)  dB(t)+\int_{%
%TCIMACRO{\U{211d} }%
%BeginExpansion
\mathbb{R}
%EndExpansion
_{0}}K^{\prime}\left(  t,t,e\right)  \tilde{N}(dt,de)\\
-\left(  \int_{t}^{T}\frac{\partial Z^{\prime}}{\partial t}\left(  t,s\right)
dB(s)\right)  dt-\left(  \int_{t}^{T}\int_{%
%TCIMACRO{\U{211d} }%
%BeginExpansion
\mathbb{R}
%EndExpansion
_{0}}\frac{\partial K^{\prime}}{\partial t}\left(  t,s,e\right)  \tilde
{N}(dt,de)\right)  dt,
\end{array}
\label{eq3.29}%
\end{equation}
where we have denoted by $\nabla$ the partial derivatives w.r.t. $x,y,z$ and
$u$ and the Fr\'{e}chet derivative w.r.t $k$ such that $\nabla=\left(
\frac{\partial}{\partial x},\frac{\partial}{\partial y},\frac{\partial
}{\partial z},\nabla_{k},\frac{\partial}{\partial u}\right)  ^{t}$ with the
second Fr\'{e}chet derivative $\nabla_{k}^{2}:=\nabla_{k}\nabla_{k}.$

\begin{theorem}
[Necessary maximum principle]Let $\hat{u}\in\mathcal{U}_{\mathbb{G}}$ with
corresponding solutions $\hat{X}(t),(\hat{Y}(t),\hat{Z}(t,s),\hat{K}\left(
t,s,\cdot\right)  ),
%$\newline$
\hat{\lambda}(t),\left(  \hat{p}(t),\hat{q}(t),\hat{r}\left(  t,\cdot\right)
\right)  $ of equations \eqref{eq3.6},\eqref{eq2.13a},\eqref{eq3.2} and
\eqref{eq3.4}, respectively. Then, the following are equivalent:

\begin{description}
\item[(i)]
\[
\frac{d}{d\epsilon}J\left(  \hat{u}+\epsilon\mu\right)  \mid_{\epsilon=0}=0,
\]
for all bounded $\mu$ of the form \eqref{eq3.1}.

\item[(ii)]
\[
\mathbb{E}\left[  \frac{\partial\mathcal{H}}{\partial u}(t)\mid\mathcal{G}%
_{t}\right]  _{u=\hat{u}} =0\text{ for all }t\in\left[  0,T\right]  .
\]

\end{description}
\end{theorem}

\noindent{Proof.} \quad Consider%

\begin{equation}
\frac{d}{d\epsilon}J\left(  \hat{u}+\epsilon\mu\right)  |\ _{\epsilon=0}%
=I_{1}+I_{2}+I_{3},\label{i0}%
\end{equation}
where%
\begin{align}
I_{1} &  =\mathbb{E}\left[  \int_{0}^{T}\left\{  \frac{\partial f}{\partial
x}(t)X^{\prime}\left(  t\right)  +\frac{\partial f}{\partial y}(t)Y^{\prime
}\left(  t\right)  +\frac{\partial f}{\partial u}(t)\mu\left(  t\right)
\right\}  dt\right]  ,\label{i1}\\
I_{2} &  =\mathbb{E}\left[  \varphi^{\prime}\left(  X\left(  T\right)
\right)  X^{\prime}\left(  T\right)  \right]  \nonumber\\
&  =\mathbb{E}\left[  p\left(  T\right)  X^{\prime}\left(  T\right)  \right]
-\mathbb{E}\left[  \lambda\left(  T\right)  \eta^{\prime}\left(  X\left(
T\right)  X^{\prime}\left(  T\right)  \right)  \right]  \text{,}\nonumber\\
I_{3} &  =\mathbb{E}\left[  \psi^{\prime}\left(  Y\left(  0\right)  \right)
Y^{\prime}\left(  0\right)  \right]  .\nonumber
\end{align}
By the It\^{o} formula
\begin{align*}
&  \mathbb{E}\left[  p\left(  T\right)  X^{\prime}\left(  T\right)  \right]
\\
&  =\mathbb{E}\left[  \int_{0}^{T}p\left(  t\right)  \left(  \frac{\partial
b}{\partial x}(t,t)X^{\prime}\left(  t\right)  +\frac{\partial b}{\partial
u}(t,t)\mu\left(  t\right)  \right)  dt\right.  \\
&  +\int_{0}^{T}p\left(  t\right)  \left\{  \int_{0}^{t}\left(  \frac
{\partial^{2}b}{\partial t\partial x}(t,s)X^{\prime}\left(  s\right)
+\frac{\partial^{2}b}{\partial t\partial u}(t,s)\mu\left(  s\right)  \right)
ds\right\}  dt\\
&  +\int_{0}^{T}p\left(  t\right)  \left\{  \int_{0}^{t}\left(  \frac
{\partial^{2}\sigma}{\partial t\partial x}(t,s)X^{\prime}\left(  s\right)
+\frac{\partial^{2}\sigma}{\partial t\partial u}(t,s)\mu\left(  s\right)
\right)  dB\left(  s\right)  \right\}  dt\\
&  +\int_{0}^{T}p\left(  t\right)  \left\{  \int_{0}^{t}\int_{%
%TCIMACRO{\U{211d} }%
%BeginExpansion
\mathbb{R}
%EndExpansion
_{0}}\left(  \frac{\partial^{2}\theta}{\partial t\partial x}\left(
t,s,e\right)  X^{\prime}\left(  s\right)  +\frac{\partial^{2}\theta}{\partial
t\partial u}\left(  t,s,e\right)  \mu(s)\right)  \tilde{N}(ds,de)\right\}
dt\\
&  -\int_{0}^{T}X^{\prime}\left(  t\right)  \frac{\partial\mathcal{H}%
}{\partial x}\left(  t\right)  dt+\int_{0}^{T}q\left(  t\right)  \left(
\frac{\partial\sigma}{\partial x}(t,t)X^{\prime}\left(  t\right)
+\frac{\partial\sigma}{\partial u}(t,t)\mu\left(  t\right)  \right)  dt\\
&  \left.  +\int_{0}^{T}\left(  \int_{%
%TCIMACRO{\U{211d} }%
%BeginExpansion
\mathbb{R}
%EndExpansion
_{0}}\left(  \frac{\partial\theta}{\partial x}\left(  t,t,e\right)  X^{\prime
}\left(  t\right)  +\frac{\partial\theta}{\partial u}\left(  t,t,e\right)
\mu(t)\right)  r\left(  t,e\right)  \nu\left(  de\right)  \right)  dt\right]
.
\end{align*}
From \eqref{eq2.11}, \eqref{eq2.12} and \eqref{eq2.17}, we have%
\begin{align}
&  \mathbb{E}\left[  p\left(  T\right)  X^{\prime}\left(  T\right)  \right]
\nonumber\\
&  =\mathbb{E}\left[  \int_{0}^{T}p\left(  t\right)  \left(  \frac{\partial
b}{\partial x}(t,t)X^{\prime}\left(  t\right)  +\frac{\partial b}{\partial
u}(t,t)\mu\left(  t\right)  \right)  dt\right.  \nonumber\\
&  +\int_{0}^{T}\int_{t}^{T}p\left(  s\right)  \left\{  \left(  \frac
{\partial^{2}b}{\partial s\partial x}(s,t)X^{\prime}\left(  t\right)
+\frac{\partial^{2}b}{\partial s\partial u}(s,t)\mu\left(  t\right)  \right)
ds\right\}  dt\nonumber\\
&  +\int_{0}^{T}\left\{  \left(  \frac{\partial^{2}\sigma}{\partial s\partial
x}(s,t)X^{\prime}\left(  t\right)  +\frac{\partial^{2}\sigma}{\partial
s\partial u}(s,t)\mu\left(  t\right)  \right)  \int_{t}^{T}\mathbb{E}\left[
D_{t}p\left(  s\right)  \mid\mathcal{F}_{t}\right]  ds\right\}  dt\nonumber\\
&  +%
%TCIMACRO{\dint _{0}^{T}}%
%BeginExpansion
{\displaystyle\int_{0}^{T}}
%EndExpansion
\left\{
%TCIMACRO{\dint \nolimits_{t}^{T}}%
%BeginExpansion
{\displaystyle\int\nolimits_{t}^{T}}
%EndExpansion%
%TCIMACRO{\dint \limits_{\mathbb{R}_{0}}}%
%BeginExpansion
{\displaystyle\int\limits_{\mathbb{R}_{0}}}
%EndExpansion
\left(  \frac{\partial^{2}\theta}{\partial s\partial x}\left(  s,t,e\right)
X^{\prime}\left(  t\right)  +\frac{\partial^{2}\theta}{\partial s\partial
u}\left(  s,t,e\right)  \mu(t)\right)  \mathbb{E}\left[  D_{t,e}p\left(
s\right)  \mid\mathcal{F}_{t}\right]  \nu\left(  de\right)  \right\}
dt\nonumber\\
&  -%
%TCIMACRO{\dint _{0}^{T}}%
%BeginExpansion
{\displaystyle\int_{0}^{T}}
%EndExpansion
\frac{\partial\mathcal{H}}{\partial x}\left(  t\right)  X^{\prime}\left(
t\right)  dt+\int_{0}^{T}\left(  \frac{\partial\sigma}{\partial x}%
(t,t)X^{\prime}\left(  t\right)  +\frac{\partial\sigma}{\partial u}%
(t,t)\mu\left(  t\right)  \right)  q\left(  t\right)  dt\nonumber\\
&  \left.  +%
%TCIMACRO{\dint _{0}^{T}}%
%BeginExpansion
{\displaystyle\int_{0}^{T}}
%EndExpansion
\left(
%TCIMACRO{\dint \limits_{\mathbb{R}_{0}}}%
%BeginExpansion
{\displaystyle\int\limits_{\mathbb{R}_{0}}}
%EndExpansion
\left(  \frac{\partial\theta}{\partial x}\left(  t,t,e\right)  X^{\prime
}\left(  t\right)  +\frac{\partial\theta}{\partial u}\left(  t,t,e\right)
\mu(t)\right)  r\left(  t,e\right)  \nu\left(  de\right)  \right)  dt\right]
.\label{i2}%
\end{align}
By the It\^{o} formula and \eqref{eq3.28}-\eqref{eq3.29}, we get%
\[%
\begin{array}
[c]{l}%
\mathbb{E}\left[  \psi^{\prime}\left(  Y\left(  0\right)  \right)  Y^{\prime
}\left(  0\right)  \right]  =\mathbb{E}\left[  \lambda\left(  0\right)
Y^{\prime}\left(  0\right)  \right]  \\
=\mathbb{E}\left[  \lambda\left(  T\right)  Y^{\prime}\left(  T\right)
\right]  +\mathbb{E}\left[
%TCIMACRO{\dint _{0}^{T}}%
%BeginExpansion
{\displaystyle\int_{0}^{T}}
%EndExpansion
\lambda\left(  t\right)  \left\{  \nabla g(t,t)\left(  X^{\prime}\left(
t\right)  ,Y^{\prime}\left(  t\right)  ,Z^{\prime}\left(  t,t\right)
,K^{\prime}\left(  t,t,\cdot\right)  ,\mu\left(  t\right)  \right)
^{t}\right.  \right.  \\
-%
%TCIMACRO{\dint _{t}^{T}}%
%BeginExpansion
{\displaystyle\int_{t}^{T}}
%EndExpansion
\left\{  \nabla\left(  \dfrac{\partial g}{\partial t}(t,s),\dfrac{\partial
g}{\partial z}(t,s)\dfrac{\partial Z}{\partial t}\left(  t,s\right)
,\nabla_{k}g(t,s)\dfrac{\partial K}{\partial t}\left(  t,s,\cdot\right)
\right)  \right.  \\
\left.  \left(  X^{\prime}\left(  t\right)  ,Y^{\prime}\left(  t\right)
,Z^{\prime}\left(  t,t\right)  ,K^{\prime}\left(  t,t,\cdot\right)
,\mu\left(  t\right)  \right)  ^{t}\right\}  dsdt\\
+%
%TCIMACRO{\dint _{0}^{T}}%
%BeginExpansion
{\displaystyle\int_{0}^{T}}
%EndExpansion
\lambda\left(  t\right)  \left(
%TCIMACRO{\dint _{t}^{T}}%
%BeginExpansion
{\displaystyle\int_{t}^{T}}
%EndExpansion
\dfrac{\partial Z^{\prime}}{\partial t}\left(  t,s\right)  dB(s)\right)  dt+%
%TCIMACRO{\dint _{0}^{T}}%
%BeginExpansion
{\displaystyle\int_{0}^{T}}
%EndExpansion
\lambda\left(  t\right)  \left(
%TCIMACRO{\dint _{t}^{T}}%
%BeginExpansion
{\displaystyle\int_{t}^{T}}
%EndExpansion%
%TCIMACRO{\dint _{\mathbb{R}}}%
%BeginExpansion
{\displaystyle\int_{\mathbb{R}}}
%EndExpansion
\dfrac{\partial K^{\prime}}{\partial t}\left(  t,s,e\right)  \tilde
{N}(dt,de)\right)  dt\\
\left.  -%
%TCIMACRO{\dint _{0}^{T}}%
%BeginExpansion
{\displaystyle\int_{0}^{T}}
%EndExpansion
\dfrac{\partial\mathcal{H}}{\partial y}\left(  t\right)  Y^{\prime}\left(
t\right)  dt-%
%TCIMACRO{\dint _{0}^{T}}%
%BeginExpansion
{\displaystyle\int_{0}^{T}}
%EndExpansion
\dfrac{\partial\mathcal{H}}{\partial z}\left(  t\right)  Z^{\prime}\left(
t,s\right)  dt\right.  \\
\left.  -%
%TCIMACRO{\dint \nolimits_{0}^{T}}%
%BeginExpansion
{\displaystyle\int\nolimits_{0}^{T}}
%EndExpansion%
%TCIMACRO{\dint _{\mathbb{R}_{0}}}%
%BeginExpansion
{\displaystyle\int_{\mathbb{R}_{0}}}
%EndExpansion
\dfrac{d\nabla_{k}\mathcal{H}}{\partial\nu}\left(  t\right)  K^{\prime}\left(
t,s,e\right)  \nu\left(  de\right)  dt\right]  .
\end{array}
\]
From \eqref{eq2.15}-\eqref{eq.2.19} and the Fubini Theorem, we have
\begin{align}
&  \mathbb{E}\left[  \psi^{\prime}\left(  Y\left(  0\right)  \right)
Y^{\prime}\left(  0\right)  \right]  =\mathbb{E}\left[  \lambda\left(
T\right)  Y^{\prime}\left(  T\right)  \right]  \nonumber\\
&  +\mathbb{E}\left[  \int_{0}^{T}\lambda\left(  t\right)  \left\{  \nabla
g(t,t)\left(  X^{\prime}\left(  t\right)  ,Y^{\prime}\left(  t\right)
,Z^{\prime}\left(  t,t\right)  ,K^{\prime}\left(  t,t,\cdot\right)
,\mu\left(  t\right)  \right)  ^{t}\right\}  \right.  \nonumber\\
&  +\int_{0}^{T}\int_{0}^{t}\lambda\left(  s\right)  \left\{  \nabla\left(
\dfrac{\partial g}{\partial t}(s,t),\dfrac{\partial g}{\partial z}%
(s,t)\dfrac{\partial Z}{\partial t}\left(  s,t\right)  ,\nabla_{k}%
g(s,t)\dfrac{\partial K}{\partial t}\left(  s,t,\cdot\right)  \right)
\right.  \nonumber\\
&  \left.  \left(  X^{\prime}\left(  s\right)  ,Y^{\prime}\left(  s\right)
,Z^{\prime}\left(  t,s\right)  ,K^{\prime}\left(  t,s,\cdot\right)
,\mu\left(  s\right)  \right)  ^{t}\right\}  dsdt\nonumber\\
&  -\int_{0}^{T}\frac{\partial\mathcal{H}}{\partial y}\left(  t\right)
Y^{\prime}\left(  t\right)  dt-\int_{0}^{T}\frac{\partial\mathcal{H}}{\partial
z}\left(  t\right)  Z^{\prime}\left(  t,s\right)  dt\nonumber\\
&  \text{ }\left.  -%
%TCIMACRO{\dint \nolimits_{0}^{T}}%
%BeginExpansion
{\displaystyle\int\nolimits_{0}^{T}}
%EndExpansion%
%TCIMACRO{\dint _{\mathbb{R}_{0}}}%
%BeginExpansion
{\displaystyle\int_{\mathbb{R}_{0}}}
%EndExpansion
\dfrac{d\nabla_{k}\mathcal{H}}{\partial\nu}\left(  t\right)  K^{\prime}\left(
t,s,e\right)  \nu\left(  de\right)  dt\right]  .\label{i3}%
\end{align}
Using that%
\begin{equation}%
\begin{array}
[c]{ll}%
\dfrac{\partial\mathcal{H}}{\partial x}\left(  t\right)   & =\dfrac{\partial
f}{\partial x}\left(  t\right)  +\frac{\partial b}{\partial x}\left(
t,t\right)  p\left(  t\right)  +\dfrac{\partial\sigma}{\partial x}\left(
t,t\right)  q\left(  t\right)  +\lambda\left(  t\right)  \dfrac{\partial
g}{\partial x}(t,t)\\
& +%
%TCIMACRO{\dint \limits_{\mathbb{R}_{0}}}%
%BeginExpansion
{\displaystyle\int\limits_{\mathbb{R}_{0}}}
%EndExpansion
\dfrac{\partial\theta}{\partial x}\left(  t,t,e\right)  r\left(  t,e\right)
\nu\left(  de\right)  +%
%TCIMACRO{\dint _{0}^{t}}%
%BeginExpansion
{\displaystyle\int_{0}^{t}}
%EndExpansion
\dfrac{\partial^{2}g}{\partial s\partial x}(s,t)\lambda\left(  s\right)  ds\\
& +%
%TCIMACRO{\dint _{t}^{T}}%
%BeginExpansion
{\displaystyle\int_{t}^{T}}
%EndExpansion
\dfrac{\partial^{2}b}{\partial s\partial x}\left(  s,t\right)  p\left(
s\right)  ds+%
%TCIMACRO{\dint _{t}^{T}}%
%BeginExpansion
{\displaystyle\int_{t}^{T}}
%EndExpansion
\dfrac{\partial^{2}\sigma}{\partial s\partial x}\left(  s,t\right)
\mathbb{E}\left[  D_{t}p\left(  s\right)  \mid\mathcal{F}_{t}\right]  ds\\
& +%
%TCIMACRO{\dint \limits_{\mathbb{R}_{0}}}%
%BeginExpansion
{\displaystyle\int\limits_{\mathbb{R}_{0}}}
%EndExpansion
\dfrac{\partial^{2}\theta}{\partial s\partial x}\left(  s,t,e\right)
\mathbb{E}\left[  D_{t,e}p\left(  s\right)  \mid\mathcal{F}_{t}\right]
\nu(de)ds\\
& +%
%TCIMACRO{\dint _{0}^{t}}%
%BeginExpansion
{\displaystyle\int_{0}^{t}}
%EndExpansion
\dfrac{\partial^{2}g}{\partial x\partial z}(s,t)\dfrac{\partial Z}{\partial
s}\left(  s,t\right)  \lambda\left(  s\right)  ds\\
& +%
%TCIMACRO{\dint _{0}^{t}}%
%BeginExpansion
{\displaystyle\int_{0}^{t}}
%EndExpansion
\dfrac{\partial}{\partial x}\left(  \nabla_{k}g(s,t)\right)  \left(
\dfrac{\partial K}{\partial s}\left(  s,t,\cdot\right)  \right)  ds,
\end{array}
\label{i4}%
\end{equation}
and that
\begin{align}
\nabla_{k}\mathcal{H}(t) &  =\nabla_{k}g(t,t)\lambda(t)+\int_{0}^{t}\nabla
_{k}\left(  \frac{\partial}{\partial s}g(s,t)\right)  \lambda(s)ds\nonumber\\
&  +\int_{0}^{t}\nabla_{k}\left(  \frac{\partial g}{\partial z}\left(
s,t\right)  \right)  \frac{\partial Z}{\partial s}\left(  s,t\right)
\lambda\left(  s\right)  ds\nonumber\\
&  +\int_{0}^{t}\nabla_{k}^{2}g\left(  s,t\right)  \frac{\partial K}{\partial
s}\left(  s,t,\cdot\right)  \lambda\left(  s\right)  ds,\label{i5}%
\end{align}
similarly for $\dfrac{\partial\mathcal{H}}{\partial y}\left(  t\right)  $ and
$\dfrac{\partial\mathcal{H}}{\partial z}\left(  t\right)  $. Combining
$\left(  \ref{i1}\right)  -\left(  \ref{i3}\right)  $\ with $\left(
\ref{i0}\right)  ,\left(  \ref{i4}\right)  -\left(  \ref{i5}\right)  $ and by
the definition of $\mu,$ we obtain%
\[
\frac{d}{d\epsilon}J\left(  u+\epsilon\mu\right)  \mid_{\epsilon=0}%
=\mathbb{E}\left[  \int_{0}^{T}\frac{\partial\mathcal{H}}{\partial u}\left(
t\right)  \mu\left(  t\right)  dt\right]  =\mathbb{E}\left[  \int
_{t}^{t+\varepsilon}\frac{\partial\mathcal{H}}{\partial u}\left(  s\right)
ds\alpha\right]  .
\]
We conclude that%
\[
\frac{d}{d\epsilon}J\left(  u+\epsilon\mu\right)  \mid_{\epsilon=0}=0
\]
if and only if%
\[
\mathbb{E}\left[  \frac{\partial\mathcal{H}}{\partial u}(t)\mid\mathcal{G}%
_{t}\right]  =0.
\]

\hfill$\square$

\section{Existence and uniqueness of solutions of BSVIE}

In order to prove existence and uniqueness solution of the backward stochastic
Volterra integral equations (BSVIE), let us introduce the following BSVIE in
the unknown $Y,Z$ and $K$:
\begin{equation}%
\begin{array}
[c]{c}%
Y(t)=\zeta\left(  t\right)  +\int_{t}^{T}g(t,s,Y(s),Z(t,s),K(t,s,\cdot
))ds-\int_{t}^{T}Z(t,s)dB(s)\\
-\int_{t}^{T}\int_{%
%TCIMACRO{\U{211d} }%
%BeginExpansion
\mathbb{R}
%EndExpansion
_{0}}K(t,s,e)\tilde{N}(ds,de),t\in\left[  0,T\right]  .
\end{array}
\label{a3}%
\end{equation}
In this section we prove existence and uniqueness of solutions of $(\ref{a3}%
)$, following the approach by Yong \cite{Y} and \cite{Yo}, but now we have
jumps. The papers by Wang and Zhang \cite{WZ},\ and by Ren \cite{Ren} studied
more general cases of $(\ref{a3})$ and our case can be seen as a particular
case of theirs, but we have included this part because it will be more
convenient for the reader to have a direct and simple approach. For related
results on BSVIE, we refer to Shi and Wang and Yong \cite{SW}-\cite{SWY}.
\newline Let us now introduce the following spaces:\newline For any $\beta
\geq0$, let $\triangle:=\left\{  \left(  t,s\right)  \in\left[  0,T\right]
^{2}:t\leq s\right\}  $ and $H_{\triangle}^{2,\beta}\left[  0,T\right]  $ be a
space of all processes $\left(  Y,Z,K\right)  $, such that $Y:\left[
0,T\right]  \times\Omega\rightarrow%
%TCIMACRO{\U{211d} }%
%BeginExpansion
\mathbb{R}
%EndExpansion
$ is $\mathbb{F}$-adapted$,$ and $Z:\triangle\times\Omega\rightarrow%
%TCIMACRO{\U{211d} }%
%BeginExpansion
\mathbb{R}
%EndExpansion
$, $K:\triangle\times%
%TCIMACRO{\U{211d} }%
%BeginExpansion
\mathbb{R}
%EndExpansion
_{0}\times\Omega\rightarrow%
%TCIMACRO{\U{211d} }%
%BeginExpansion
\mathbb{R}
%EndExpansion
$ with $s\mapsto Z(t,s)$ and $s\mapsto K\left(  t,s,\cdot\right)  $ being
$\mathbb{F}$-adapted on $[t,T],$ equipped with the norm
\begin{align*}
&  \left\Vert \left(  Y,Z,K\right)  \right\Vert _{H_{\triangle}^{2,\beta
}\left[  0,T\right]  }^{2}\\
&  :=\mathbb{E}\int\nolimits_{0}^{T}\left[  e^{\beta t}\left\vert Y\left(
t\right)  \right\vert ^{2}+\int\nolimits_{t}^{T}e^{\beta s}\left\vert Z\left(
t,s\right)  \right\vert ^{2}ds+\int\nolimits_{t}^{T}\int_{%
%TCIMACRO{\U{211d} }%
%BeginExpansion
\mathbb{R}
%EndExpansion
_{0}}e^{\beta s}\left\vert K(t,s,e)\right\vert ^{2}\nu(ds,de)\right]  dt.
\end{align*}
Clearly $H_{\triangle}^{2,\beta}\left[  0,T\right]  $ is a Hilbert space. It
is easy to see that for any $\beta>0$ , the norm $\left\Vert .\right\Vert
_{H_{\triangle}^{2,\beta}\left[  0,T\right]  }$ is equivalent to $\left\Vert
.\right\Vert _{H_{\triangle}^{2,0}\left[  0,T\right]  }$ obtained from
$\left\Vert .\right\Vert _{H_{\triangle}^{2,\beta}\left[  0,T\right]  }$ by
taking $\beta$ $=0$. We now make the following assumptions:

\textbf{Assumptions (H.1)}

\begin{itemize}
\item The function $g:\left[  0,T\right]  ^{2}\times%
%TCIMACRO{\U{211d} }%
%BeginExpansion
\mathbb{R}
%EndExpansion
^{3}\times L^{2}(\nu)\times\Omega\rightarrow%
%TCIMACRO{\U{211d} }%
%BeginExpansion
\mathbb{R}
%EndExpansion
$, is such that
\end{itemize}

\begin{enumerate}
\item $\mathbb{E}\left[  \int\limits_{0}^{T}\left(  \int\limits_{t}%
^{T}g(t,s,0,0,0)ds\right)  ^{2}dt\right]  <+\infty,$

\item There exists a constant $c>0$, such that, for all $t,s\in\left[
0,T\right]  $%
\[%
\begin{array}
[c]{l}%
\left\vert g(t,s,y,z,k(\cdot))-g(t,s,y^{\prime},z^{\prime},k^{\prime}%
(\cdot))\right\vert \\
\leq c\left(  \left\vert y-y^{\prime}\right\vert +\left\vert z-z^{\prime
}\right\vert +\left(  \int_{%
%TCIMACRO{\U{211d} }%
%BeginExpansion
\mathbb{R}
%EndExpansion
_{0}}\left\vert k(e)-k^{\prime}\left(  e\right)  \right\vert ^{2}%
\nu(de)\right)  ^{\frac{1}{2}}\right)
\end{array}
\]
for all $y,y^{\prime},z,z^{\prime},k(\cdot),k^{\prime}(\cdot).$
\end{enumerate}

\begin{itemize}
\item $\zeta\left(  \cdot\right)  \in L_{\mathcal{F}_{T}}^{2}\left(  \Omega,%
%TCIMACRO{\U{211d} }%
%BeginExpansion
\mathbb{R}
%EndExpansion
\right)  .$
\end{itemize}

\begin{theorem}
Under the assumptions (H.1), there exists a unique solution $(Y,Z,K)\in
H_{\triangle}^{2,\beta}\left[  0,T\right]  $ of the BSVIE $(\ref{a3}).$
\end{theorem}

For a given triple of processes $\left(  y\left(  \cdot\right)  ,z\left(
\cdot,\cdot\right)  ,k\left(  \cdot,\cdot,\cdot\right)  \right)  \in
H_{\triangle}^{2,\beta}\left[  0,T\right]  ,$ consider the following simple
BSVIE in the unknown triple $(Y,Z,K)$:
\begin{equation}
Y\left(  t\right)  =\zeta\left(  t\right)  +\int_{t}^{T}\bar{g}\left(
t,s\right)  ds-\int_{t}^{T}Z(t,s)dB(s)-\int_{t}^{T}\int_{%
%TCIMACRO{\U{211d} }%
%BeginExpansion
\mathbb{R}
%EndExpansion
_{0}}K(t,s,e)\tilde{N}(ds,de), \label{a4}%
\end{equation}
where we denote by%
\[
\bar{g}\left(  t,s\right)  =g\left(  t,s,y\left(  s\right)  ,z\left(
t,s\right)  ,k\left(  t,s,\cdot\right)  \right)  ,\text{ for }\left(
t,s\right)  \in\triangle.
\]
To solve $(\ref{a4})$ for $(Y,Z,K)$, we introduce the following family of BSDE
(parameterized by $t\in\lbrack0,T])$:%
\[
\bar{Y}\left(  r,t\right)  =\zeta\left(  t\right)  +\int_{r}^{T}\bar{g}\left(
t,s\right)  ds-\int_{r}^{T}\bar{Z}(s,t)dB(s)-\int_{r}^{T}\int_{%
%TCIMACRO{\U{211d} }%
%BeginExpansion
\mathbb{R}
%EndExpansion
_{0}}\bar{K}(s,t,e)\tilde{N}(ds,de),\text{ }r\in(t,T],
\]
It is well-known that the above BSDE admits a unique adapted solution $\left(
\bar{Y}\left(  \cdot,t\right)  ,\bar{Z}(\cdot,t),\bar{K}(\cdot,t,\cdot
)\right)  $ and the following estimate holds:%
\[%
\begin{array}
[c]{l}%
\mathbb{E}\left[  \underset{r\in\left[  t,T\right]  }{\sup}\left\vert \bar
{Y}\left(  r,t\right)  \right\vert ^{2}+\int_{t}^{T}\left\vert \bar
{Z}(s,t)\right\vert ^{2}ds+\int\nolimits_{t}^{T}\int_{%
%TCIMACRO{\U{211d} }%
%BeginExpansion
\mathbb{R}
%EndExpansion
_{0}}\left\vert \bar{K}(s,t,e)\right\vert ^{2}\nu(de)ds\right] \\
\leq C\mathbb{E}\left[  \left\vert \zeta\left(  t\right)  \right\vert
^{2}+\left(  \int_{t}^{T}\bar{g}\left(  t,s\right)  ds\right)  ^{2}\right]  .
\end{array}
\]
Now let%
\[
Y\left(  t\right)  =\bar{Y}\left(  t,t\right)  ,\text{ }Z\left(  t,s\right)
=\bar{Z}(s,t),\text{ }K\left(  t,s,\cdot\right)  =\bar{K}(s,t,\cdot),\text{
for all }\left(  t,s\right)  \in\triangle.
\]
Then $(Y(\cdot),Z(\cdot,\cdot),K\left(  \cdot,\cdot\right)  )$ is an adapted
solution to the BSVIE $(\ref{a4})$, and%
\begin{align*}
&  \mathbb{E}\left[  \left\vert Y\left(  t\right)  \right\vert ^{2}+\int
_{t}^{T}\left\vert Z(t,s)\right\vert ^{2}ds+\int\nolimits_{t}^{T}\int_{%
%TCIMACRO{\U{211d} }%
%BeginExpansion
\mathbb{R}
%EndExpansion
_{0}}\left\vert K(t,s,e)\right\vert ^{2}\nu(de)ds\right] \\
&  =\mathbb{E}\left[  \left\vert \zeta\left(  t\right)  +\int_{t}^{T}\bar
{g}\left(  t,s\right)  ds\right\vert ^{2}\right] \\
&  \leq2\mathbb{E}\left[  \left\vert \zeta\left(  t\right)  \right\vert
^{2}+\left(  \int_{t}^{T}\bar{g}\left(  t,s\right)  ds\right)  ^{2}\right]  .
\end{align*}
Therefore, by integrating both sides of the inequality above, we get%
\begin{align*}
&  \mathbb{E}\left[  \int_{0}^{T}\left(  \left\vert Y\left(  t\right)
\right\vert ^{2}+\int_{t}^{T}\left\vert Z(t,s)\right\vert ^{2}ds+\int
\nolimits_{t}^{T}\int_{%
%TCIMACRO{\U{211d} }%
%BeginExpansion
\mathbb{R}
%EndExpansion
_{0}}\left\vert K(t,s,e)\right\vert ^{2}\nu(de)ds\right)  dt\right] \\
&  \leq2\mathbb{E}\int_{0}^{T}\left[  \left\vert \zeta\left(  t\right)
\right\vert ^{2}+\left(  \int_{t}^{T}\bar{g}\left(  t,s\right)  ds\right)
^{2}\right]  dt.
\end{align*}
Adding and subtracting $g\left(  t,s,0,0,0\right)  $ on the left side, then by
the Lipschitz assumption, we obtain
\begin{align*}
&  \mathbb{E}\left[  \int_{0}^{T}\left(  \left\vert Y\left(  t\right)
\right\vert ^{2}+\int_{t}^{T}\left\vert Z(t,s)\right\vert ^{2}ds+\int
\nolimits_{t}^{T}\int_{%
%TCIMACRO{\U{211d} }%
%BeginExpansion
\mathbb{R}
%EndExpansion
_{0}}\left\vert K(t,s,e)\right\vert ^{2}\nu(de)ds\right)  dt\right] \\
&  \leq C\mathbb{E}\int\nolimits_{0}^{T}\left[  \left\vert \zeta\left(
t\right)  \right\vert ^{2}+\left(  \int\nolimits_{t}^{T}g\left(
t,s,0,0,0\right)  ds\right)  ^{2}\right]  dt\\
&  +C\mathbb{E}\left[  \int\nolimits_{0}^{T}\left(  \left\vert y\left(
t\right)  \right\vert ^{2}+\int\nolimits_{t}^{T}\left\vert z(s)\right\vert
^{2}ds+\int\nolimits_{t}^{T}\int_{%
%TCIMACRO{\U{211d} }%
%BeginExpansion
\mathbb{R}
%EndExpansion
_{0}}\left\vert k(t,s,e)\right\vert ^{2}\nu(de)ds\right)  dt\right]  ,
\end{align*}
for some constant $C$. Thus, $(y,z,k)\mapsto(Y,Z,K)$ defines a map from
$H_{\triangle}^{2,\beta}\left[  0,T\right]  $ to itself. \newline Now, we want
to prove that this mapping is contracting in $H_{\triangle}^{2,\beta}\left[
0,T\right]  $ under the norm $\left\Vert .\right\Vert _{H_{\triangle}%
^{2,\beta}\left[  0,T\right]  }.$ We show that if for $i=1,2$, $(y_{i}%
,z_{i},k_{i})\in H_{\triangle}^{2,\beta}\left[  0,T\right]  $ and
$(Y_{i},Z_{i},K_{i})$ is the corresponding adapted solution to equation
$(\ref{a3})$, then%
\[%
\begin{array}
[c]{l}%
\mathbb{E}\left[  \int_{0}^{T}\left(  e^{\beta t}\left\vert Y_{1}\left(
t\right)  -Y_{2}\left(  t\right)  \right\vert ^{2}+\int_{t}^{T}e^{\beta
s}\left\vert Z_{1}(t,s)-Z_{2}(t,s)\right\vert ^{2}ds\right.  \right. \\
\left.  \left.  +\int\nolimits_{t}^{T}e^{\beta s}\int_{%
%TCIMACRO{\U{211d} }%
%BeginExpansion
\mathbb{R}
%EndExpansion
_{0}}\left\vert K_{1}(t,s,e)-K_{2}(t,s,e)\right\vert ^{2}\nu(de)ds\right)
dt\right] \\
\leq\dfrac{C}{\beta}\mathbb{E}\left[  \int\nolimits_{0}^{T}\left(  e^{\beta
t}\left\vert y_{1}\left(  t\right)  -y_{2}\left(  t\right)  \right\vert
^{2}+\int\nolimits_{t}^{T}e^{\beta s}\left\vert z_{1}(t,s)-z_{2}%
(t,s)\right\vert ^{2}ds\right.  \right. \\
\left.  \left.  +\int\nolimits_{t}^{T}e^{\beta s}\int_{%
%TCIMACRO{\U{211d} }%
%BeginExpansion
\mathbb{R}
%EndExpansion
_{0}}\left\vert k_{1}(t,s,e)-k_{2}(t,s,e)\right\vert ^{2}\nu(de)ds\right)
dt\right]  ,
\end{array}
\]
which means that%
\[
\left\Vert \left(  Y,Z,K\right)  \right\Vert _{H_{\triangle}^{2,\beta}\left[
0,T\right]  }^{2}\leq\dfrac{C}{\beta}\left\Vert \left(  y,z,k\right)
\right\Vert _{H_{\triangle}^{2,\beta}\left[  0,T\right]  }^{2}.
\]
Hence, the mapping $(y,z,k)\mapsto(Y,Z,K)$ is contracting on $H_{\triangle
}^{2,\beta}\left[  0,T\right]  $ for large enough $\beta>0$. Then, $(Y,Z,K)$
is a unique solution for the BSVIE $(\ref{a3})$.\qquad$\qquad\qquad
\qquad\qquad\square$ \bigskip

\section{Application: Optimal recursive utility consumption}

As an illustration of our general results above, we now apply them to solve
the optimal recursive utility consumption problem \eqref{eq1.5} described in
the Introduction. Our example is related to the examples discussed in
\cite{AO2} and \cite{OS}, but now the cash flow is modelled by a stochastic
Volterra equation and the utility is represented by the recursive utility. As
pointed out after \eqref{eq1.2} in the Introduction, the Volterra equation
contains history terms and can therefore be viewed as a model for a system
with memory. Thus, we assume that the cash flow $X(t)=X^{c}(t)$ being exposed
to a $\mathbb{G}$-adapted consumption rate $c(t)$, satisfies the stochastic
Volterra equation%

\begin{equation}%
\begin{array}
[c]{c}%
X(t)=\xi+\int_{0}^{t}\left(  \alpha(t,s)-c(s)\right)  X(s)ds+\int_{0}^{t}%
\beta(t,s)X(s)dB(s)\\
\text{ \ \ \ \ \ \ \ }+\int_{0}^{t}\int_{%
%TCIMACRO{\U{211d} }%
%BeginExpansion
\mathbb{R}
%EndExpansion
_{0}}\pi(t,s,e)X(s)\tilde{N}(ds,de),t\in\left[  0,T\right]  ,
\end{array}
\label{vs}%
\end{equation}
where we assume for simplicity that $\xi$ is a (deterministic) constant and
$\alpha$, $\beta:\left[  0,T\right]  ^{2}\rightarrow%
%TCIMACRO{\U{211d} }%
%BeginExpansion
\mathbb{R}
%EndExpansion
$ and $\pi:\left[  0,T\right]  ^{2}\times%
%TCIMACRO{\U{211d} }%
%BeginExpansion
\mathbb{R}
%EndExpansion
_{0}\rightarrow%
%TCIMACRO{\U{211d} }%
%BeginExpansion
\mathbb{R}
%EndExpansion
$ are deterministic functions with $\alpha$, $\beta$ and $\pi$
bounded.$\smallskip$The FSVIE $\left(  \ref{vs}\right)  $ can be written in
its differential form as%

\[%
\begin{array}
[c]{l}%
dX(t)=\left(  \alpha(t,t)-c(t)\right)  X(t)dt+\left(  \int_{0}^{t}%
\frac{\partial\alpha}{\partial t}(t,s)X(s)ds\right)  dt\\
+\beta(t,t)X(t)dB(t)+\left(  \int_{0}^{t}\frac{\partial\beta}{\partial
t}(t,s)X(s)dB(s)\right)  dt\\
+\int_{%
%TCIMACRO{\U{211d} }%
%BeginExpansion
\mathbb{R}
%EndExpansion
_{0}}\pi(t,t,e)X(t)\tilde{N}(dt,de)+\left(  \int_{%
%TCIMACRO{\U{211d} }%
%BeginExpansion
\mathbb{R}
%EndExpansion
_{0}}\int_{0}^{t}\frac{\partial\pi}{\partial t}(t,s,e)X(s)\tilde
{N}(ds,de)\right)  dt,t\in\left[  0,T\right]  .
\end{array}
\]
The recursive utility process $Y(t)$ of Duffie and Epstein \cite{DE} has the
following linear form%

\begin{equation}%
\begin{array}
[c]{c}%
dY(t)=-\left[  \gamma(t)Y(t)+\ln c(t)X(t)\right]  dt+Z(t)dB(t)\\
\text{ \ \ \ \ \ \ \ \ \ \ \ \ \ \ \ }+\int_{%
%TCIMACRO{\U{211d} }%
%BeginExpansion
\mathbb{R}
%EndExpansion
_{0}}K(t,e)\tilde{N}(dt,de),t\in\left[  0,T\right]  .
\end{array}
\label{bs}%
\end{equation}
Our problem \eqref{eq1.5} is to maximise the performance functional
\[
J(c):=Y^{c}(0)
\]
over all control processes $c\in\mathcal{U}_{\mathbb{G}}$, where in this case
$\mathcal{U}_{\mathbb{G}}$ is the set of all $\mathbb{G}$-adapted nonnegative
processes. \newline This problem is a special case of the problem discussed in
the previous sections, with $f=0$, $\varphi=0$ and $\psi(y)=y.$ The
Hamiltonian associated to our problem is defined by%

\[%
\begin{array}
[c]{l}%
\mathcal{H}(t,s,x,y,p,q)=\left(  \alpha(t,t)-c(t)\right)  px+\int_{t}^{T}%
\frac{\partial\alpha}{\partial s}(s,t)x(s)p(s)ds\\
+\beta(t,t)qx+\int_{t}^{T}\frac{\partial\beta}{\partial s}(s,t)x(s)\mathbb{E}%
\left[  \left.  D_{t}p(s)\right\vert \mathcal{F}_{t}\right]  ds\\
+\int_{%
%TCIMACRO{\U{211d} }%
%BeginExpansion
\mathbb{R}
%EndExpansion
_{0}}\pi(t,t,e)xr(t,e)\nu(de)\\
+\int_{%
%TCIMACRO{\U{211d} }%
%BeginExpansion
\mathbb{R}
%EndExpansion
_{0}}\int_{t}^{T}\frac{\partial\pi}{\partial s}(s,t,e)x(s)\mathbb{E}\left[
\left.  D_{t,e}p(s)\right\vert \mathcal{F}_{t}\right]  \nu(de)ds\\
+\left[  \gamma(t)y+\ln c(t)+\ln x\right]  \lambda.
\end{array}
\]
The corresponding backward-forward system for the adjoint processes $(p,q,r)$
and $\lambda$ are%

\begin{equation}
\left\{
\begin{array}
[c]{l}%
dp(t)=-\left[  \left(  \alpha(t,t)-c(t)\right)  p(t)+\int_{t}^{T}%
\frac{\partial\alpha}{\partial s}(s,t)p(s)ds\right. \\
+\beta(t,t)q(t)+\int_{t}^{T}\frac{\partial\beta}{\partial s}(s,t)\mathbb{E}%
\left[  \left.  D_{t}p(s)\right\vert \mathcal{F}_{t}\right]  ds+\int_{%
%TCIMACRO{\U{211d} }%
%BeginExpansion
\mathbb{R}
%EndExpansion
_{0}}\pi(t,t,e)r(t,e)\nu(de)\\
\left.  +\int_{%
%TCIMACRO{\U{211d} }%
%BeginExpansion
\mathbb{R}
%EndExpansion
_{0}}\int_{t}^{T}\frac{\partial\pi}{\partial s}(s,t,e)\mathbb{E}\left[
\left.  D_{t,e}p(s)\right\vert \mathcal{F}_{t}\right]  \nu(ds,de)+\frac
{\lambda(t)}{X(t)}\right]  dt\\
+q(t)dB(t)+\int_{%
%TCIMACRO{\U{211d} }%
%BeginExpansion
\mathbb{R}
%EndExpansion
_{0}}r(t,e)\tilde{N}(dt,de),t\in\left[  0,T\right]  ,\\
p(T)=0,
\end{array}
\right.  \label{b}%
\end{equation}
and%
\begin{equation}
\left\{
\begin{array}
[c]{l}%
d\lambda(t)=\gamma(t)\lambda(t)dt,t\in\left[  0,T\right]  ,\\
\lambda(0)=1.
\end{array}
\right.  \label{eq4.5}%
\end{equation}
The solution of the differential equation \eqref{eq4.5} is
\[
\lambda(t)=\exp\left(  -\int_{0}^{t}\gamma(s)ds\right)  ,t\in\left[
0,T\right]  .
\]
Now, maximising the Hamiltonian w.r.t $c$ gives the first order condition%

\begin{equation}
c(t)=\mathbb{E}\left[  \left.  \frac{\lambda(t)}{p(t)X(t)}\right\vert
\mathcal{G}_{t}\right]  ,\quad t\in\lbrack0,T].\label{eq4.6a}%
\end{equation}
Applying It\^{o}'s formula, we get%
\[%
\begin{array}
[c]{l}%
d\left(  p(t)X(t)\right)  =p(t)dX(t)+X(t)dp(t)+d\left[  p(t)X(t)\right]  \\
=p(t)\left\{  \left(  \alpha(t,t)-c(t)\right)  X(t)dt+\left(  \int_{0}%
^{t}\frac{\partial\alpha}{\partial t}(t,s)X(s)ds\right)  dt\right.  \\
+\beta(t,t)X(t)dB(t)+\left(  \int_{0}^{t}\frac{\partial\beta}{\partial
t}(t,s)X(s)dB(s)\right)  dt\\
\left.  +\int_{%
%TCIMACRO{\U{211d} }%
%BeginExpansion
\mathbb{R}
%EndExpansion
_{0}}\pi(t,t,e)X(t)\tilde{N}(dt,de)+\left(  \int_{%
%TCIMACRO{\U{211d} }%
%BeginExpansion
\mathbb{R}
%EndExpansion
_{0}}\int_{0}^{t}\frac{\partial\pi}{\partial t}(t,s,e)X(s)\tilde
{N}(ds,de)\right)  dt\right\}  \\
-X(t)\left\{  \left(  \alpha(t,t)-c(t)\right)  p(t)dt+\left(  \int_{t}%
^{T}\frac{\partial\alpha}{\partial s}(s,t)p(s)ds\right)  dt+\beta
(t,t)q(t)dt\right.  \\
+\left(  \int_{t}^{T}\frac{\partial\beta}{\partial s}(s,t)\mathbb{E}\left[
\left.  D_{t}p(s)\right\vert \mathcal{F}_{t}\right]  ds\right)  dt+\int_{%
%TCIMACRO{\U{211d} }%
%BeginExpansion
\mathbb{R}
%EndExpansion
_{0}}\pi(t,t,e)r(t,e)\nu(dt,de)\\
+\left(  \int_{%
%TCIMACRO{\U{211d} }%
%BeginExpansion
\mathbb{R}
%EndExpansion
_{0}}\int_{t}^{T}\frac{\partial\pi}{\partial s}(s,t,e)\mathbb{E}\left[
\left.  D_{t,e}p(s)\right\vert \mathcal{F}_{t}\right]  \nu(de)ds\right)
dt+\frac{\lambda(t)}{X(t)}dt\\
\left.  +q(t)dB(t)+\int_{%
%TCIMACRO{\U{211d} }%
%BeginExpansion
\mathbb{R}
%EndExpansion
_{0}}r(t,e)\tilde{N}(dt,de)\right\}  \\
+\beta(t,t)X(t)q(t)dt+\int_{%
%TCIMACRO{\U{211d} }%
%BeginExpansion
\mathbb{R}
%EndExpansion
_{0}}\pi(t,t,e)X(t)r(t,e)\nu(dt,de).
\end{array}
\]
Collecting the terms, we see that the above reduces to
\[
\left\{
\begin{array}
[c]{l}%
p(t)X(t)=p(0)X(0)-\int_{0}^{t}\lambda(s)ds\\
+\int_{0}^{t}\left\{  p(s)X(s)\beta(s,s)-X(s)q(s)\right\}  dB(s)\\
+\int_{0}^{t}\int_{%
%TCIMACRO{\U{211d} }%
%BeginExpansion
\mathbb{R}
%EndExpansion
_{0}}\left\{  p(s)X(s)\pi(s,s,e)-X(s)r(s,e)\right\}  \tilde{N}(ds,de),t\in
\left[  0,T\right]  ,\\
p(T)X(T)=0.
\end{array}
\right.
\]
Therefore, if we define
\[%
\begin{array}
[c]{l}%
P(t)=p(t)X(t),\\
Q(t)=p(s)X(s)\beta(s,s)-X(s)q(s),\\
R(t,e)=p(s)X(s)\pi(s,s,e)-X(s)r(s,e),
\end{array}
\]
then $(P,Q,R)$ solves the linear BSDE%
\[
\left\{
\begin{array}
[c]{l}%
dP(t)=-\lambda(t)dt+Q(t)dB(t)+\int_{%
%TCIMACRO{\U{211d} }%
%BeginExpansion
\mathbb{R}
%EndExpansion
_{0}}R(t,e)\tilde{N}(dt,de),t\in\left[  0,T\right]  ,\\
P(T)=0.
\end{array}
\right.
\]
The solution of this linear BSDE is%
\[
P(t)=\mathbb{E}\left[  \int_{t}^{T}\left.  \lambda(s)ds\right\vert
\mathcal{F}_{t}\right]  =p(t)X(t).
\]
Combined with \eqref{eq4.6a} this gives
\begin{equation}
c(t)=c^{\ast}(t)=\mathbb{E}\left[  \left.  \frac{\exp\left(  -\int_{0}%
^{t}\gamma(s)ds\right)  }{\mathbb{E}\left[  \int_{t}^{T}\left.  \exp\left(
-\int_{0}^{s}\gamma(r)dr\right)  ds\right\vert \mathcal{F}_{t}\right]
}\right\vert \mathcal{G}_{t}\right]  .\label{eq4.6}%
\end{equation}
In particular, since $\lambda>0$ by \eqref{eq4.5} we get that $p(t)X(t)>0$.
Thus we see that $c(t)$ is well-defined in \eqref{eq4.6a} and $c^{\ast}(t)>0$
for all $t\in\lbrack0,T]$. Therefore $c^{\ast}\in\mathcal{U}_{\mathbb{G}}$,
and we conclude that $c^{\ast}$ is indeed optimal. We have proved

\begin{theorem}
The optimal recursive utility consumption rate $c^{*}(t)$ for the problem
\eqref{eq1.5} (with $\xi$ constant) is given by \eqref{eq4.6}.
\end{theorem}

\section{Appendix}

\subsection{Some basic concepts from Banach space theory}

%{ Throughout this paper, \textcolor{blue}{we  assume that
%the probability space $\Omega$  %(introduced in Section 2)
%is a Banach space}, {where $(\Omega,   \FC, \mathbb{F}:=\{\FC_t\}_{t \geq 0}, P)$ is as introduced in the beginning of Section 1.
%\ref{motiv}.
%Throughout this paper
%$\mathbb{G}  :=\{\mathcal{G}_t\}_{t \geq 0}$
%is a given subfiltration of $\mathbb{F}$.}
%}
To explain the notation used in this paper, we briefly recall some basic
concepts from Banach space theory:\newline

Let $\mathcal{X},\mathcal{Y}$ be two Banach spaces with norms $\|
\cdot\|_{\mathcal{X}},\| \cdot\|_{\mathcal{Y}}$, respectively, and let $F :
\mathcal{X} \rightarrow\mathcal{Y}$.

\begin{itemize}
\item We say that $F$ has a directional derivative (or G\^ateaux derivative)
at $v\in\mathcal{X}$ in the direction $w\in\mathcal{X}$ if
\[
D_{w}F(v):=\lim_{\varepsilon\rightarrow0}\frac{1}{\varepsilon}(F(v+\varepsilon
w)-F(v))
\]
exists.

\item We say that $F$ is Fr\' echet differentiable at $v \in\mathcal{X}$ if
there exists a continuous linear map $A: \mathcal{X} \rightarrow\mathcal{Y}$
such that
\[
\lim_{\substack{h \rightarrow0 \\h \in\mathcal{X}}} \frac{1}{
\|h\|_{\mathcal{X}}} \| F(v+h) - F(v) - A(h)\|_{\mathcal{Y}} = 0.
\]
In this case we call $A$ the \textit{gradient} (or Fr\' echet derivative) of
$F$ at $v$ and we write
\[
A =\nabla_{v} F.
\]

\item If $F$ is Fr\' echet differentiable at $v$ with Fr\' echet derivative
$\nabla_{v}F$, then $F$ has a directional derivative in all directions
$w\in\mathcal{X}$ and
\[
D_{w}F(v):=\left\langle \nabla_{v}F,w\right\rangle =\nabla_{v}F(w)=\nabla
_{v}Fw.
\]

\end{itemize}

In particular, note that if $F$ is a linear operator, then $\nabla_{v}F=F$ for
all $v$.\newline

\subsection{A brief review of Hida-Malliavin calculus for L\'{e}vy processes}

For the convenience of the reader, in this section we recall the basic
definition and properties of Hida-Malliavin calculus for L\'{e}vy processes
related to this paper. The following summary is based on \cite{AO}. A general
reference for this presentation is the book \cite{DOP}.\newline First, recall
the L\'{e}vy--It\^{o} decomposition theorem, which states that any L\'{e}vy
process $Y(t)$ with
\[
\mathbb{E}[Y^{2}(t)]<\infty\quad\mbox{for all}\quad t
\]
can be written
\[
Y(t)=at+bB(t)+\int_{0}^{t}\int_{\mathbb{R}_{0}}e\tilde{N}(ds,de)
\]
with constants $a$ and $b$. In view of this we see that it suffices to deal
with Hida-Malliavin calculus for $B(\cdot)$ and for
\[
\chi(\cdot):=\int_{0}^{\cdot}\int_{\mathbb{R}_{0}}e\tilde{N}(ds,de)
\]
separately.

\subsection{Hida-Malliavin calculus for $B(\cdot)$}

A natural starting point is the Wiener-It\^{o} chaos expansion theorem, which
states that any $F\in L^{2}(\mathcal{F}_{T},P)$ can be written
\begin{equation}
F=\sum_{n=0}^{\infty}I_{n}(f_{n})\label{a}%
\end{equation}
for a unique sequence of symmetric deterministic functions $f_{n}\in
L^{2}(\rho^{n})$, where $\rho$ is Lebesgue measure on $[0,T]$ and
\[
I_{n}(f_{n})=n!\int_{0}^{T}\int_{0}^{t_{n}}\cdots\int_{0}^{t_{2}}f_{n}%
(t_{1},\cdots,t_{n})dB(t_{1})dB(t_{2})\cdots dB(t_{n})
\]
(the $n$-times iterated integral of $f_{n}$ with respect to $B(\cdot)$) for
$n=1,2,\ldots$ and $I_{0}(f_{0})=f_{0}$ when $f_{0}$ is a constant.\newline
Moreover, we have the isometry
\[
\mathbb{E}[F^{2}]=||F||_{L^{2}(P)}^{2}=\sum_{n=0}^{\infty}n!||f_{n}%
||_{L^{2}(\rho^{n})}^{2}.
\]

\begin{definition}
[Hida-Malliavin derivative $D_{t}$ with respect to $B(\cdot)$]
%{Definition 2.1}
\hfill\break\textrm{Let $\mathbb{D}_{1,2}^{(B)}$ be the space of all $F\in
L^{2}(\mathcal{F}_{T},P)$ such that its chaos expansion \eqref{a} satisfies
\[
||F||_{\mathbb{D}_{1,2}^{(B)}}^{2}:=\sum_{n=1}^{\infty}nn!||f_{n}%
||_{L^{2}(\rho^{n})}^{2}<\infty.
\]
}

\textrm{For $F\in\mathbb{D}^{(B)}_{1,2}$ and $t\in[0,T]$, we define the
\emph{Hida-Malliavin derivative} or \emph{the stochastic gradient}) of $F$ at
$t$ (with respect to $B(\cdot)$), $D_{t}F,$ by
\begin{align}
\label{eq2.5a}D_{t}F=\sum^{\infty}_{n=1}nI_{n-1}(f_{n}(\cdot,t)),
\end{align}
where the notation $I_{n-1}(f_{n}(\cdot,t))$ means that we apply the
$(n-1)$-times iterated integral to the first $n-1$ variables $t_{1},\cdots,
t_{n-1}$ of $f_{n}(t_{1},t_{2},\cdots,t_{n})$ and keep the last variable
$t_{n}=t$ as a parameter.}
\end{definition}

One can easily check that
\[
\mathbb{E}\Big[\int_{0}^{T}(D_{t}F)^{2}dt\Big]=\sum_{n=1}^{\infty}%
nn!||f_{n}||_{L^{2}(\rho^{n})}^{2}=||F||_{\mathbb{D}_{1,2}^{(B)}}^{2},
\]
so $(t,\omega)\longmapsto D_{t}F(\omega)$ belongs to $L^{2}(\rho\times P)$.

%\underline{Example 2.2}

\begin{example}
\textrm{If $F=\int_{0}^{T}f(t)dB(t)$ with $f\in L^{2}(\rho)$ deterministic,
then
\[
D_{t}F=f(t)\mbox{ for }a.a.\,t\in\lbrack0,T].
\]
More generally, if }$\mathrm{\Psi}$\textrm{$(s)$ is Skorohod integrable,
}$\mathrm{\Psi}$\textrm{$(s)\in\mathbb{D}_{1,2}$ for $a.a.\;s$ and
$D_{t}\mathrm{\Psi}(s)$ is Skorohod integrable for $a.a.\;t$, then
\[
D_{t}\Big(\int_{0}^{T}\mathrm{\Psi}(s)\delta B(s)\Big)=\int_{0}^{T}%
D_{t}\mathrm{\Psi}(s)\delta B(s)+\mathrm{\Psi}%
(t)\;\mbox{for a.a. $(t,\omega)$},
\]
where $\int_{0}^{T}\mathrm{\Psi}(s)\delta B(s)$ denotes the Skorohod integral
of a process }$\mathrm{\Psi}$\textrm{ with respect to $B(\cdot)$. }
\end{example}

Some other basic properties of the Hida-Malliavin derivative $D_{t}$ are the following:

\begin{enumerate}
\item[(i)] \textbf{Chain rule } \newline Suppose $F_{1},\ldots,F_{m}%
\in\mathbb{D}_{1,2}^{(B)}$ and that $\mathrm{\Psi}:\mathbb{R}^{m}%
\rightarrow\mathbb{R}$ is $C^{1}$
%$e^\bot$
with bounded partial derivatives. Then, $\mathrm{\Psi}(F_{1},\cdots,F_{m}%
)\in\mathbb{D}_{1,2}$ and
\[
D_{t}\mathrm{\Psi}(F_{1},\cdots,F_{m})=\sum_{i=1}^{m}\frac{\partial
\mathrm{\Psi}}{\partial x_{i}}(F_{1},\cdots,F_{m})D_{t}F_{i}.
\]

\item[(ii)] \textbf{Duality formula} \newline Suppose $\mathrm{\Psi}(t)$ is
$\mathcal{F}_{t}-$adapted with $\mathbb{E}[\int_{0}^{T}\mathrm{\Psi}%
^{2}(t)dt]<\infty$ and let $F\in\mathbb{D}_{1,2}^{(B)}$. Then,
\begin{equation}
\mathbb{E}[F\int_{0}^{T}\mathrm{\Psi}(t)dB(t)]=\mathbb{E}[\int_{0}%
^{T}\mathrm{\Psi}(t)D_{t}Fdt].\label{n}%
\end{equation}

\item[(iii)] \textbf{Malliavin derivative and adapted processes}\newline If
$\mathrm{\Psi}$ is an $\mathbb{F}$-adapted process, then%
\[
D_{s}\mathrm{\Psi}(t)=0\text{ for }s>t.
\]

\end{enumerate}

\begin{remark}
We put $D_{t}\mathrm{\Psi}(t)=\underset{s\rightarrow t-}{\lim}D_{s}%
\mathrm{\Psi}(t)$ (if the limit exists).
\end{remark}

\begin{remark}
It was proved in \cite{AaOPU} that one can extend the Hida-Malliavin
derivative operator $D_{t}$ from $\mathbb{D}_{1,2}$ to all of $L^{2}%
(\mathcal{F}_{T},P)$ in such a way that, also denoting the extended operator
by $D_{t}$, for all $F\in L^{2}(\mathcal{F}_{T},P)$ we have
\begin{equation}
D_{t}F\in(\mathcal{S})^{\ast}\text{ and }(t,\omega)\mapsto\mathbb{E}%
[D_{t}F\mid\mathcal{F}_{t}]\text{ belongs to }L^{2}(\rho\times
P)\label{eq2.10a}%
\end{equation}
Here $(\mathcal{S})^{\ast}$ is the Hida space of stochastic
distributions.\newline Moreover, the following \emph{generalized
Clark-Haussmann-Ocone formula} was proved:
\begin{equation}
F=\mathbb{E}[F]+\int_{0}^{T}\mathbb{E}[D_{t}F\mid\mathcal{F}_{t}%
]dB(t)\label{eq2.11a}%
\end{equation}
for all $F\in L^{2}(\mathcal{F}_{T},P)$. See Theorem 3.11 in \cite{AaOPU} and
also Theorem 6.35 in \cite{DOP}.\newline We can use this to get the following
extension of the duality formula \eqref{n}:
\end{remark}

\begin{proposition}
\textbf{The generalized duality formula}\newline Let $F\in L^{2}%
(\mathcal{F}_{T},P)$ and let $\mathrm{\Psi}(t,\omega)\in L^{2}(\rho\times P)$
be adapted. Then
\[
\mathbb{E}[F\int_{0}^{T}\mathrm{\Psi}(t)dB(t)]=\mathbb{E}[\int_{0}%
^{T}\mathbb{E}[D_{t}F\mid\mathcal{F}_{t}]\mathrm{\Psi}(t)dt].
\]

\end{proposition}

{\noindent\textit{Proof\quad}} By \eqref{eq2.10a} and \eqref{eq2.11a} and the
It\^{o} isometry we get
\begin{align*}
&  \mathbb{E}[F\int_{0}^{T}\mathrm{\Psi}(t)dB(t)]=\mathbb{E}[(\mathbb{E}%
[F]+\int_{0}^{T}\mathbb{E}[D_{t}F\mid\mathcal{F}_{t}]dB(t))(\int_{0}%
^{T}\mathrm{\Psi}(t)dB(t))]\\
&  =\mathbb{E}[\int_{0}^{T}\mathbb{E}[D_{t}F\mid\mathcal{F}_{t}]\mathrm{\Psi
}(t)dt].
\end{align*}
\hfill$\square$

\subsection{Hida-Malliavin calculus for $\tilde N(\cdot)$}

The construction of a stochastic derivative/Hida-Malliavin derivative in the
pure jump martingale case follows the same lines as in the Brownian motion
case. In this case, the corresponding Wiener-It\^{o} Chaos Expansion Theorem
states that any $F\in L^{2}(\mathcal{F}_{T},P)$ (where, in this case,
$\mathcal{F}_{t}=\mathcal{F}_{t}^{(\tilde{N})}$ is the $\sigma-$algebra
generated by $\chi(s):=\int_{0}^{s}\int_{\mathbb{R}_{0}}e\tilde{N}%
(dr,de);\;0\leq s\leq t$) can be written as
\begin{equation}
F=\sum_{n=0}^{\infty}I_{n}(f_{n});\;f_{n}\in\hat{L}^{2}((\rho\times\nu
)^{n}),\label{c}%
\end{equation}
where $\hat{L}^{2}((\rho\times\nu)^{n})$ is the space of functions
$f_{n}(t_{1},e_{1},\ldots,t_{n},e_{n})$, $t_{i}\in\lbrack0,T]$, $e_{i}%
\in\mathbb{R}_{0}$ such that $f_{n}\in L^{2}((\rho\times\nu)^{n})$ and $f_{n}$
is symmetric with respect to the pairs of variables $(t_{1},\rho_{1}%
),\ldots,(t_{n},\rho_{n}).$\newline It is important to note that in this case,
the $n-$times iterated integral $I_{n}(f_{n})$ is taken with respect to
$\tilde{N}(dt,de)$ and not with respect to $d\chi(t).$ Thus, we define
\[
I_{n}(f_{n}):=n!\int_{0}^{T}\!\!\int_{\mathbb{R}_{0}}\!\int_{0}^{t_{n}}%
\!\int_{\mathbb{R}_{0}}\cdots\int_{0}^{t_{2}}\!\!\int_{\mathbb{R}_{0}}%
f_{n}(t_{1},e_{1},\cdots,t_{n},e_{n})\tilde{N}(dt_{1},de_{1})\cdots\tilde
{N}(dt_{n},de_{n})
\]
for $f_{n}\in\hat{L}^{2}((\rho\times\nu)^{n}).$\newline The It\^{o} isometry
for stochastic integrals with respect to $\tilde{N}(dt,de)$ then gives the
following isometry for the chaos expansion:
\[
||F||_{L^{2}(P)}^{2}=\sum_{n=0}^{\infty}n!||f_{n}||_{L^{2}((\rho\times\nu
)^{n})}^{2}.
\]
As in the Brownian motion case, we use the chaos expansion to define the
Malliavin derivative. Note that in this case, there are two parameters $t,e,$
where $t$ represents time and $e\neq0$ represents a generic jump size.

\begin{definition}
[Hida-Malliavin derivative $D_{t,e}$ with respect to $\tilde{N}(\cdot,\cdot)$%
]Let $\mathbb{D}_{1,2}^{(\tilde{N})}$ be the space of all $F\in L^{2}%
(\mathcal{F}_{T},P)$ such that its chaos expansion \eqref{c} satisfies
\[
||F||_{\mathbb{D}_{1,2}^{(\tilde{N})}}^{2}:=\sum_{n=1}^{\infty}nn!||f_{n}%
||_{L^{2}((\rho\times\nu)^{2})}^{2}<\infty.
\]
For $F\in\mathbb{D}_{1,2}^{(\tilde{N})}$, we define the Hida-Malliavin
derivative of $F$ at $(t,e)$ (with respect to $\tilde{N}(\cdot))$, $D_{t,e}F,$
by
\[
D_{t,e}F:=\sum_{n=1}^{\infty}nI_{n-1}(f_{n}(\cdot,t,e)),
\]
where $I_{n-1}(f_{n}(\cdot,t,e))$ means that we perform the $(n-1)-$times
iterated integral with respect to $\tilde{N}$ to the first $n-1$ variable
pairs $(t_{1},e_{1}),\cdots,(t_{n},e_{n}),$ keeping $(t_{n},e_{n})=(t,e)$ as a
parameter.
%}

\end{definition}

In this case, we get the isometry.
\[
\mathbb{E}[\int_{0}^{T}\int_{\mathbb{R}_{0}}(D_{t,e}F)^{2}\nu(de)dt]=\sum
_{n=0}^{\infty}nn!||f_{n}||_{L^{2}((\rho\times\nu)^{n})}^{2}=||F||_{\mathbb{D}%
_{1,2}^{(\tilde{N})}}^{2}.
\]

\begin{example}
\textrm{If $F=\int_{0}^{T}\int_{\mathbb{R}_{0}}f(t,e)\tilde{N}(dt,de)$ for
some deterministic $f(t,e)\in L^{2}(\rho\times\nu)$, then
\[
D_{t,e}F=f(t,e)\mbox{ for }a.a.\,(t,e).
\]
More generally, if }$\mathrm{\Phi}$\textrm{$(s,e)$ is Skorohod integrable with
respect to $\tilde{N}(\delta s,de)$, }$\mathrm{\Phi}$\textrm{$(s,e)\in
\mathbb{D}_{1,2}^{(\tilde{N})}$ for $a.a.\,s,e$ and $D_{t,z}\Phi(s,e)$ is
Skorohod integrable for $a.a.\,(t,z)$, then
\[
D_{t,z}(\int_{0}^{T}\!\int_{\mathbb{R}_{0}}\mathrm{\Phi(s,e)}\tilde{N}(\delta
s,de))=\int_{0}^{T}\int_{\mathbb{R}_{0}}D_{t,z}\mathrm{\Phi(s,e)}\tilde
{N}(\delta s,de)+\mathrm{\Phi}(t,z)\;\mbox{ for }a.a.\,t,z,
\]
where $\int_{0}^{T}\int_{\mathbb{R}_{0}}\mathrm{\Phi(s,e)}\tilde{N}(\delta
s,de)$ denotes the \emph{Skorohod integral\/} of }$\mathrm{\Phi}$\textrm{ with
respect to $\tilde{N}(\cdot,\cdot).$ (See \cite{DOP} for a definition of such
Skorohod integrals and for more details.) }
\end{example}

The properties of $D_{t,e}$ corresponding to those of $D_{t}$ are the following:

\begin{itemize}
\item[(i)] \textbf{Chain rule}\newline Suppose $F_{1},\cdots,F_{m}%
\in\mathbb{D}_{1,2}^{(\tilde{N})}$ and that $\phi:\mathbb{R}^{m}%
\rightarrow\mathbb{R}$ is continuous and bounded. Then, $\phi(F_{1}%
,\cdots,F_{m})\in\mathbb{D}_{1,2}^{(\tilde{N})}$ and
\[
D_{t,e}\phi(F_{1},\cdots,F_{m})=\phi(F_{1}+D_{t,e}F_{1},\ldots,F_{m}%
+D_{t,e}F_{m})-\phi(F_{1},\ldots,F_{m}).
\]

\item[(ii)] \textbf{Duality formula }

Suppose $\mathrm{\Phi}(t,e)$ is $\mathcal{F}_{t}$-adapted and $\mathbb{E}%
[\int_{0}^{T}\int_{\mathbb{R}_{0}}\mathrm{\Phi}^{2}(t,e)\nu(de)dt]<\infty$ and
let $F\in\mathbb{D}_{1,2}^{(\tilde{N})}$. Then,
\[
\mathbb{E}\Big[F\int_{0}^{T}\int_{\mathbb{R}_{0}}\mathrm{\Phi}(t,e)\tilde
{N}(dt,de)\Big]=\mathbb{E}\Big[\int_{0}^{T}\int_{\mathbb{R}_{0}}\mathrm{\Phi
}(t,e)D_{t,e}F\nu(de)dt\Big].
\]

\item[(iii)] \textbf{Hida-Malliavin derivative and adapted processes}\newline
If $\mathrm{\Phi}$ is an $\mathbb{F}$-adapted process, then,
\[
D_{s,e}\mathrm{\Phi}(t)=0\text{ for all }s>t.
\]

\end{itemize}

\begin{remark}
We put $D_{t,e}\mathrm{\Phi}(t)=\underset{s\rightarrow t-}{\lim}%
D_{s,e}\mathrm{\Phi}(t)$ ( if the limit exists).
\end{remark}

\begin{remark}
As in Remark 3.2 we note that there is an extension of the Hida-Malliavin
derivative $D_{t,e}$ from $\mathbb{D}_{1,2}^{(\tilde{N})}$ to $L^{2}%
(\mathcal{F}_{t}\times P)$ such that the following extension of the duality
theorem holds:

\begin{proposition}
\textbf{Generalized duality formula}\newline Suppose $\mathrm{\Phi}(t,e)$ is
$\mathcal{F}_{t}$-adapted and $\mathbb{E}[\int_{0}^{T}\int_{\mathbb{R}_{0}%
}\mathrm{\Phi}^{2}(t,e)\nu(de)dt]<\infty$ and let $F\in L^{2}(\mathcal{F}%
_{T}\times P)$. Then,
\[
\mathbb{E}\Big[F\int_{0}^{T}\int_{\mathbb{R}_{0}}\mathrm{\Phi}(t,e)\tilde
{N}(dt,de)\Big]=\mathbb{E}\Big[\int_{0}^{T}\int_{\mathbb{R}_{0}}\mathrm{\Phi
}(t,e)\mathbb{E}[D_{t,e}F\mid\mathcal{F}_{t}]\nu(de)dt\Big].
\]

\end{proposition}

We refer to Theorem 13.26 in \cite{DOP}.\newline\ \newline We emphasize that
this generalized Hida-Malliavin derivative $DX$ exists for all $X\in L^{2}(P)$
as an element of the Hida stochastic distribution space $(\mathcal{S})^{\ast}%
$, and it has the property that the conditional expectation $\mathbb{E}%
[DX|\mathcal{F}_{t}]$ belongs to $L^{2}(\rho\times P)$, where $\rho$ is
Lebesgue measure on $[0,T]$. Therefore, when using this generalized
Hida-Malliavin derivative, combined with conditional expectation, no
assumptions on Hida-Malliavin differentiability in the classical sense are
needed; we can work on the whole space of random variables in $L^{2}(P)$.
\end{remark}

\section{Acknowledgments}

We want to thank Yaozhong Hu and Yanqing Wang helpful comments. We are also
grateful to an anonymous referee for a very valuable and comprehensive report,
which helped us to improve the paper considerably.


\begin{thebibliography}{99}                                                                                               %


\bibitem {AaOPU}Aase, K., \O ksendal, B., Privault, N. and Ub\o e, J.: White
noise generalizations of the Clark-Haussmann-Ocone theorem with application to
mathematical finance. Finance Stochast. 4 (2000), 465-496.

\bibitem {AO}Agram, N., \O ksendal, B.: Malliavin calculus and optimal control
of stochastic Volterra equations. J. Optim. Theory Appl. , DOI
10.1007/s10957-015-0735-5 (2015).

\bibitem {AO2}Agram, N., \O ksendal, B.: Infinite horizon optimal control of
forward-backward stochastic differential equations with delay. J. Comput.
Appl. Math. 259, 336--349 (2014).

\bibitem {DOP}Di Nunno, G., \O ksendal, B. and Proske, F.: Malliavin Calculus
for L\'{e}vy Processes and Applications to Finance. Corrected, Second
Printing, Springer (2009).

\bibitem {DE}Duffie, D., Epstein, L.G.: Stochastic differential utility,
Econometrica 60 (2) 353--394 (1992).

\bibitem {HO}Hu, Y. and \O ksendal, B.: Linear backward stochastic Volterra
equations. Manuscript August (2016).

\bibitem {L}Lin, J.: Adapted solution of a backward stochastic nonlinear
Volterra integral equation, Stochastic Analysis and Applications, vol. 20, no.
1, pp. 165--183 (2002).

\bibitem {os}\O ksendal, B. and Sulem, A.: Risk minimization in financial
markets modeled by It\^{ }o-L\'{ }evy processes. Afrika Matematika (2014),
DOI: 10.1007/s13370-014-02489-9.

\bibitem {OS}\O ksendal, B. and Sulem, A.: Optimal control of predictive
mean-field equations and applications to finance. In F.E. Benth and G. Di
Nunno (eds.), Stochastics of Environmental and Financial Economics, Springer
Proceedings in Mathematics and Statistics 138, DOI
10.1007/978-3-319-23425-0\_12 (2016).

\bibitem {Ren}Ren, Y.: On solutions of backward stochastic Volterra integral
equations with jumps in Hilbert spaces. J Optim Theory Appl 144: 319--333 (2010).

\bibitem {SWY}Shi, Y., Wang, T. and Yong, J.: Optimal control problems of
forward-backward stochastic Volterra integral equations. arXiv: 1404.7577v1 (2014).

\bibitem {SW}Shi, Y. and Wang, T.: Solvability of general backward stochastic
Volterra integral equations. J. Korean Math. Soc. 49 No. 6, pp. 1301--1321 (2012).

\bibitem {WS}Wang, T. and Shi, Y.: A maximum principle for forward-backward
stochastic Volterra integral equations and applications in finance. arXiv
1004.2206v1 (2010).

\bibitem {WZ1}Wang, Z. and Zhang, X.: Optimal control problems of
forward-backward stochastic Volterra integral equations with closed control
regions. arXiv:1602.05661v1 (2016).

\bibitem {WZ}Wang, Z. and Zhang, X.: Non-Lipschitz backward stochastic
Volterra type equations with jumps. Stoch. Dyn.07:479-496 (2007).

\bibitem {WX}Wei, Q. and Xiao, X.: An optimal control problem of
forward-backward stochastic Volterra integral equations with state
constraints. Abstract and Applied Analysis, Volume (2014), Article ID 432718,
16 (2014).

\bibitem {Y}Yong, Y.: Backward stochastic Volterra equations and some related
problems, Stochastic Processes and their Applications, 116, 779-795 (2006).

\bibitem {Yo}Yong, J.: Backward stochastic Volterra integral equations- a
brief survey. Appl. Math. J. Chinese Univ. 28(4): 383-394 (2013).
\end{thebibliography}
\end{document}